\documentclass[11pt]{amsart}
\usepackage{amssymb}

\begin{document}

\newtheorem{theorem}{Theorem}    
\newtheorem{proposition}[theorem]{Proposition}
\newtheorem{conjecture}[theorem]{Conjecture}
\def\theconjecture{\unskip}
\newtheorem{corollary}[theorem]{Corollary}
\newtheorem{lemma}[theorem]{Lemma}
\newtheorem{observation}[theorem]{Observation}
\theoremstyle{definition}
\newtheorem{definition}{Definition}
\newtheorem{remark}{Remark}
\def\theremark{\unskip}
\newtheorem{question}{Question}
\def\thequestion{\unskip}
\newtheorem{example}{Example}
\def\theexample{\unskip}
\newtheorem{problem}{Problem}

\numberwithin{theorem}{section}
\numberwithin{definition}{section}
\numberwithin{equation}{section}

\def\intprod{\mathbin{\lr54}}
\def\reals{{\mathbb R}}
\def\integers{{\mathbb Z}}
\def\complex{{\mathbb C}\/}
\def\naturals{{\mathbb N}\/}
\def\distance{\operatorname{distance}\,}
\def\degree{\operatorname{degree}\,}
\def\dim{\operatorname{dimension}\,}
\def\Span{\operatorname{span}\,}
\def\ZZ{ {\mathbb Z} }
\def\e{\varepsilon}
\def\p{\partial}
\def\rp{{ ^{-1} }}
\def\Re{\operatorname{Re\,} }
\def\Im{\operatorname{Im\,} }
\def\ov{\overline}
\def\bx{{\bf{x}}}
\def\eps{\varepsilon}
\def\lt{L^2}

\def\scriptx{{\mathcal X}}
\def\scriptb{{\mathcal B}}
\def\scripta{{\mathcal A}}
\def\scriptk{{\mathcal K}}
\def\scriptd{{\mathcal D}}
\def\scriptp{{\mathcal P}}
\def\scriptl{{\mathcal L}}
\def\scriptv{{\mathcal V}}
\def\scripti{{\mathcal I}}
\def\scripth{{\mathcal H}}
\def\scriptm{{\mathcal M}}
\def\scripte{{\mathcal E}}
\def\scriptt{{\mathcal T}}
\def\scriptb{{\mathcal B}}
\def\frakg{{\mathfrak g}}
\def\frakG{{\mathfrak G}}

\author{Michael Christ}
\address{
        Michael Christ\\
        Department of Mathematics\\
        University of California \\
        Berkeley, CA 94720-3840, USA}
\email{mchrist@math.berkeley.edu}
\thanks{The first author was supported in part by NSF grant DMS-9970660}

\author{Xiaochun Li}
\address{
        Xiaochun Li\\
        Department of Mathematics\\
        UCLA \\
        Los Angeles, CA 90095-1555, USA}
\email{xcli@math.ucla.edu}
\thanks{The second author was supported by NSF grant DMS-0140376}

\author{Terence Tao}
\address{
        Terence Tao\\
        Department of Mathematics\\
        UCLA \\
        Los Angeles, CA 90095-1555, USA}
\email{tao@math.ucla.edu}
\thanks{The third author is a Clay Prize Fellow and is supported by a
grant from the Packard foundation.}

\author{Christoph Thiele}
\address{
        Christoph Thiele\\
        Department of Mathematics\\
        UCLA \\
        Los Angeles, CA 90095-1555, USA}
\email{thiele@math.ucla.edu}
\thanks{The fourth author was supported by NSF grants
DMS 9970469 and DMS 9985572.}

\date{October 30, 2003}

\title[Multilinear oscillatory integrals]{On multilinear oscillatory integrals, \\nonsingular and singular}

\begin{abstract}
Basic questions concerning nonsingular multilinear operators with oscillatory
factors are posed and partially answered. 
$L^p$ norm inequalities are established for multilinear integral
operators of Calder\'on-Zygmund type which incorporate
oscillatory factors $e^{iP}$, where $P$ is a real-valued polynomial.
A related problem concerning upper bounds for measures of sublevel sets is solved. 
\end{abstract}
\maketitle

\section{Introduction}


Consider multilinear functionals 
\begin{equation} \label{defn}
\Lambda_\lambda(f_1,f_2,\cdots,f_n)
= \int_{\reals^m} e^{i\lambda P(x)} \prod_{j=1}^n f_j(\pi_j(x))\eta(x)\,dx
\end{equation}
where $\lambda\in\reals$ is a parameter, $P:\reals^m\to\reals$ is a measurable real-valued
function, $m\ge 2$, and $\eta\in C^1_0(\reals^m)$ is compactly supported. 
Each $\pi_j$ denotes the orthogonal projection from $\reals^m$
to a linear subspace $V_j\subset\reals^m$ of any dimension $\kappa\le m-1$,
and $f_j:V_j\to\complex$ is always assumed to be locally integrable
with respect to Lebesgue measure on $V_j$.
We assume for simplicity that $\kappa$ is independent of $j$. 

The integral $\Lambda_\lambda$ is well-defined if all $f_j$ belong to $L^\infty$,
and satisfies $|\Lambda_\lambda(f_1,\cdots,f_n)|\le C\prod_j\|f_j\|_{L^\infty}$.
\begin{definition} \label{defn:powerdecay}
A measurable real-valued function $P$ is said to have the power decay property,
relative to a collection of subspaces $\{V_j\}$ of $\reals^m$, 
on an open set\footnote{Our theorems
concern only polynomial phases $P$, for which we will show that nondegeneracy is
independent of $U$.}
$U\subset\reals^m$ if for any $\eta\in C^1_0(U)$ 
there exist $\eps>0$ and $C<\infty$ such that 
\begin{equation} \label{decay}
|\Lambda_\lambda(f_1,\cdots,f_n)|\le C(1+|\lambda|)^{-\epsilon}
\prod_{j=1}^n \|f_j\|_{L^\infty}
\ \text{ for all $f_j\in L^\infty(V_j)$ 
and all $\lambda\in\reals$.}
\end{equation}
Moreover, for $\eta$ supported in any fixed compact set,
$\eps$ may be taken to be independent of $\eta$, and
$C=O(\|\eta\|_{C^1})$. 
\end{definition}
The goal of this paper is to characterize those data $(P,\{V_j\})$ for which the power decay
property holds.
A necessary condition is that $P$ can not be expressed as a linear combination
of measurable functions $p_j\circ\pi_j$.
Indeed, if $P = \sum_j p_j\circ\pi_j$ then for $f_j = e^{-i\lambda p_j}\in L^\infty$
one has $e^{i\lambda P(x)}\prod_j f_j(\pi_j(x))\equiv 1$,
and there is consequently no decay.
\begin{definition}
A polynomial $P$ is said to be degenerate (relative to $\{V_j\}$)
if there exist polynomials $p_j:V_j\to\reals$
such that $P=\sum_{j=1}^n p_j\circ\pi_j$.
Otherwise $P$ is nondegenerate.
\end{definition}
\noindent
In the case $n=0$, where the collection of subspaces $\{V_j\}$ is empty,
$P$ is considered to be nondegenerate if and only if it is nonconstant.

We will show in Lemma~\ref{lemma:polyordistribution}
that degeneracy is also equivalent to the existence, in some nonempty open set,
of a representation of the more general form 
$P=\sum_{j=1}^n h_j\circ\pi_j$ where the $h_j$ are arbitrary  distributions rather than 
polynomials.

More generally, a (measurable) function $P$ is said to be degenerate if it can be expressed as
$\sum_j p_j\circ\pi_j$ for some (measurable) functions $p_j$,
or distributions. But we will restrict attention to polynomials henceforth.  

\begin{question}
Is the power decay property equivalent to nondegeneracy, for real-valued polynomials? 
\end{question}
We were led to this question by a problem concerning multilinear singular integral
operators; see Theorems~\ref{thm1} and \ref{thm2} below.
We have been able answer it, always in the affirmative, only in special cases.

Numerous variants can be formulated. 
One can ask when \eqref{decay} holds with the right-hand side replaced
by $\Theta(\lambda)\prod_{j=1}^n \|f_j\|_{L^\infty}$ 
where $\Theta(\lambda)$ is a specific function of $(1+|\lambda|)$, 
or whether there merely exists some function
$\Theta$ tending to zero as $|\lambda|\to\infty$ for which it holds. 
When \eqref{decay} does hold, one can ask what is the optimal
power $\eps$ in \eqref{decay}.  We will focus on 
the first formulation, which is the one most directly relevant
to our applications to multilinear singular integral operators,
and is possibly the most fundamental.

The extreme formulation in which all $f_j$ are measured in the strongest Lebesgue
norm $L^\infty$ is the essence of the matter.
Since $|\Lambda_\lambda(f_1,\cdots,f_n)|\le C\|f_k\|_{L^1}\prod_{j\ne k}\|f_j\|_{L^\infty}$
uniformly in $\lambda$ for any $k$, a simple interpolation argument shows that
if \eqref{decay} does hold, then a decay estimate of the same type holds
with $\prod_j \|f_j\|_\infty$ replaced by $\prod_j\|f_j\|_{p_j}$
for various $n$-tuples of indices $p_j\in(1,\infty]$. More precisely,
if the integral converges absolutely whenever each $f_j\in L^{p_j}$,
then $|\Lambda_\lambda(f_1,\cdots,f_n)|\le C(1+|\lambda|)^{-\eps}\prod_j\|f\|_{q_j}$
whenever each $q_j>p_j$, where $\eps>0$ depends on $\{p_j\},\{q_j\}$,
and other $n$-tuples $q$ can in general be obtained via further interpolations.
On the other hand, the combinations of exponents $p_j$ for which the integral
is guaranteed to converge absolutely depend strongly on $\{V_j\}$. 

When $n=0$, that is, when $\{V_j\}$ is empty,
one is dealing with oscillatory integrals 
$\Lambda_\lambda = \int e^{i\lambda P(x)}\eta(x)\,dx$,
which are of the first type in the terminology of Stein \cite{steinbook}.
The general case can also be regarded as concerning oscillatory integrals
of the first type. Indeed,
\eqref{decay} with $L^\infty$ bounds on the functions $f_j$
is equivalent to 
$|\int e^{i\lambda \phi(x)}\eta(x)\,dx| \le C|\lambda|^{-\epsilon}$
{\em uniformly} for all phase functions of the form $\phi = P-\sum_j h_j\circ\pi_j$,
where the $h_j$ are arbitrary real-valued measurable functions; one implication
is tautologous, while the second  is nearly trivial and is proved in \cite{ccw}.

The regularity condition $\eta\in C^1$ is rather arbitrary. If $\eta$ is merely
H\"older continuous then for any $s<\infty$, $\eta$ may be decomposed
as a smooth function whose $C^s$ norm is $O(|\lambda|^{C\delta})$
plus a remainder which is $O(|\lambda|^{-\delta})$ in supremum norm. If \eqref{decay}
holds for all $\eta\in C_0^s$ with a constant $C$ which
is $O(\|\eta\|_{C^s})$ then it follows from this decomposition, with
$\delta = \eps/2C$, that it continues to hold for all H\"older continuous
$\eta$.

$\prod_{j=1}^n f_j\circ\pi_j$ could be replaced
by $\prod_{j=1}^n g_j\circ\ell_j$ in \eqref{defn}, where  
each $\ell:\reals^m\to\reals^\kappa$ is linear and has full rank, 
without any increase in generality.
For any function $g\circ\ell$ may be expressed as $f\circ\pi$, 
where $\pi$ is the orthogonal projection of $\reals^m$
onto the orthocomplement $V$ of the nullspace of $\ell$; 
$f\in L^p(V)$ if and only if $g\in L^p(\reals^\kappa)$ with comparable norms. 

The linear prototype for power decay is the inequality
(see Stein \cite{steinbook}, p.\ 416 and Phong and Stein \cite{phongstein}) 
\begin{equation} \label{oscillatoryprototype}
\big|
\int e^{i\lambda P(x,y)}f(x)g(y)\eta(x,y)\,dx\,dy
\big|
=O\big(|\lambda|^{-\eps}\|f\|_{\lt}\|g\|_{\lt}\big),
\end{equation}
with suitable uniformity,
whenever $P$ is a real-valued polynomial, of bounded degree, which is nondegenerate in the sense
that some mixed partial derivative $\p^{\alpha+\beta} P/\p x^\alpha \p y^\beta$,
with both $\alpha,\beta\ne 0$, does not vanish identically;  this is equivalent to
the impossibility of decomposing $P(x,y)$ as  $p(x)+q(y)$.

A related result is the lemma
of van der Corput, which asserts for instance that $\int_a^b e^{i\lambda\varphi(t)}\,dt
= O(1+|\lambda|)^{-\eps}$ provided that $\varphi$ is real-valued and
that some derivative of $\varphi$ is bounded
away from zero (in the case of the first derivative, it is also assumed to be monotone).
Here the upper bound for the integral depends only on a lower bound for the derivative;
no upper bound on $\varphi$ or its derivatives, hence in some sense no {\it a priori}
smoothness condition, are imposed.
Higher-dimensional versions of van der Corput's lemma, likewise without assumptions
of upper bounds or extra smoothness, were established in \cite{ccw}. 

Our nondegeneracy condition replaces the hypothesis of a nonvanishing derivative. 
\eqref{decay}, like \eqref{oscillatoryprototype}, 
is invariant under the replacement of $P$ by $P-\sum_j p_j\circ\pi_j$,
for arbitrary real-valued measurable functions $p_j$, 
and there is no {\it a priori} smoothness condition or upper bound on $P-\sum_j p_j\circ\pi_j$.
Formally the power decay property more closely resembles \eqref{oscillatoryprototype}, but 
our analysis is closely related to \cite{ccw}, and combinatorial issues lurk 
here as they did there.

\section{Results} \label{section:results}
\subsection{Further definitions}
Given a degree $d$, fix a norm $\|\cdot\|_{\scriptp_d}$ 
on the finite-dimensional vector space $\scriptp_d$ of all polynomials 
in $\reals^m$ of degree $\le d$. 
Given $d$,
the norm of $P$ relative to $\{V_j\}$ is defined to be
$\inf \|P-\sum_j p_j\circ\pi_j\|_{\scriptp_d}$, where the infimum is taken
over all real-valued polynomials $p_j$ of degrees $\le d$.
Polynomials $P_\alpha$
are said to be uniformly nondegenerate if they are all of degrees $\le d$
for some finite $d$,
and if there exists $c>0$ such that the norm of $P_\alpha$
relative to $\{V_j\}$ (and relative to this degree $d$) 
is $\ge c$ for all $\alpha$.

More generally, if we are given
given a collection of subspaces $\{V_j^\alpha: \alpha\in A\}$  and a
polynomial $P_\alpha$ of uniformly bounded degree, 
we may still define uniform nondegeneracy by requiring that 
$\inf_\alpha \inf_{\{p_j\}} \|P_\alpha-\sum_j p_j\circ\pi_j^\alpha\|_{\scriptp_d}\ge c>0$,
where $\pi_j^\alpha$ is the orthogonal projection onto $V_j^\alpha$.

These definitions are independent of $d$ provided only that
each $P$ has degree $\le d$, a simple consequence of
the equivalence of any two norms on a finite-dimensional
vector space, and the fact that if $P=\sum_j p_j\circ\pi_j$
then $P=\sum_j \tilde p_j\circ\pi_j$ where $\tilde p_j$
is the sum of all terms of degrees $\le \degree(P)$
in the decomposition of $\pi_j$ as a linear combination 
of monomials.

Because the family of polynomials of degree $\le d$
has finite dimension, it is easily verified that the infimum defining the relative norm is
actually assumed for some polynomials $p_j$.
Thus either $P$ is degenerate, or the infimum is strictly positive.

\begin{definition}
A collection of subspaces $\{V_j\}$ is said to have the power
decay property if every real-valued polynomial
$P$ which is nondegenerate relative to $\{V_j\}$
has the power decay property \eqref{decay}, in every open set $U$.

A collection of subspaces $\{V_j\}$ is said to have the uniform power
decay property if \eqref{decay} holds, with uniform constants $C,\eps$,  
for any family of real-valued polynomials of bounded degrees
which are uniformly nondegenerate relative to $\{V_j\}$.
\end{definition}

A related concept turns out to be somewhat easier to analyze.
To any distribution $h\in \scriptd'(V_j)$ 
is naturally associated a distribution $h\circ\pi_j\in \scriptd'(\reals^m)$.
Denote by $V_j^\perp$ the orthogonal complement of $V_j$. 
Then $h\circ\pi_j$ is annihilated by any 
first-order constant-coefficient
differential operators of the form $w\cdot\nabla$, where $w\in V_j^\perp$.
\begin{definition}
A polynomial $P$ is said to be simply nondegenerate if there exists
a differential operator $L$ of the form $L = \prod_{j=1}^r (w_j\circ\nabla)$,
with each $w_j\in V_j^\perp$, such that $L(P)$ does not vanish identically.
\end{definition}

Simple nondegeneracy implies nondegeneracy. Indeed, for any distribution
$f$ defined on some $V_j$, the $j$-th factor of $L$ annihilates $f\circ\pi_j$,
and hence $L$ does so since all factors commute.
Thus $L(P-\sum_j f_j\circ\pi_j)\equiv L(P)$.

The converse is not true in general. Consider the case $\kappa=1$ where 
each subspace $V_j$ has dimension one. 
Let $\scriptl$ be any homogeneous constant-coefficient
linear partial differential operator with real coefficients, and let $\sigma$ be its symbol. 
Identify the dual space of $\reals^m$ with $\reals^m$ via the inner product.
Let $e$ be any unit vector and let $\pi(x) = \langle x,e\rangle e$.
Then $\scriptl$ annihilates every function $f\circ\pi$ if and only if $\sigma(e)=0$. 
To create examples of nondegenerate polynomials, fix any such nonelliptic $\scriptl$, 
choose $\{e_j\}$ to be any finite collection of distinct unit vectors satisfying $\sigma(e_j)=0$, 
and choose any real-valued polynomial $P$ such that $\scriptl(P)$ does not vanish
identically. This forces $P$ to be nondegenerate relative to $\{V_j\}$,
where $V_j$ is the span of $e_j$.
If $m\ge 3$ then $P,\scriptl,\{V_j\}$ can be chosen so that $P$ has degree two, yet
$\{V_j\}$ has arbitrarily large finite cardinality.

\subsection{Decay for nonsingular oscillatory multilinear functionals}
\label{subsect:decay}

We now state several theorems asserting that nondegeneracy implies
the power decay property, under various auxiliary hypotheses.
To formulate the first of these,
let $\{V_j: 1\le j\le n\}$  be
a collection of subspaces 
of $\reals^m$ of dimension $\kappa$. We say that they are in general position if  any subcollection
of cardinality $k\ge 1$ spans a subspace of dimension $\min(k\kappa,m)$.
It is elementary that for subspaces in general position, if each $f_j$ belongs
to $L^2(\reals)$, then the product of any $k$ functions $f_j\circ\pi_j$ belongs
to $L^2(\reals^m)$ provided that $k\kappa\le m$. 
When $k\kappa\le 2m$ the product belongs to $L^1$, by Cauchy-Schwarz.
Therefore the integral defining $\Lambda_\lambda$ converges
absolutely, and $|\Lambda_\lambda(f_1,\cdots,f_n)|\le C\prod_j\|f_j\|_{L^2}$,
uniformly in $\lambda\in\reals$.

\begin{theorem} \label{thm:uniformslices} \label{thm:slices}
Suppose that $n<2m$. Then any family $\{V_j: 1\le j\le n\}$ 
of one-dimensional subspaces of $\reals^m$ which lie in general
position has the uniform power decay property.
Moreover under these hypotheses
\begin{equation} \label{Ltwodecay}
|\Lambda_\lambda(f_1,\cdots,f_n)|
\le C(1+|\lambda|)^{-\eps}\prod_{j=1}^n\|f_j\|_{L^2}
\end{equation}
for all polynomials $P$ of bounded degree which are uniformly nondegenerate
with respect to $\{V_j\}$,
for all functions $f_j \in L^2(\reals^1)$, 
with uniform  constants $C,\eps\in\reals^+$.
\end{theorem}


The case $n=m$ is known; Phong, Stein, and Sturm \cite{phongsteinsturm} have obtained
much more precise results on the exponent in the power decay estimate, phrased 
in terms of the reduced Newton polyhedron of $P$.
The case $n<m$ follows from $n=m$. Indeed, 
choose coordinates
$(x',x'')\in \reals^{n}\times\reals^{m-n}$, in such a way that 
the first factor $\reals^{n}$ equals the span of the subspaces $V_j$.
If $\partial P/\partial x''$ does not vanish identically then \eqref{decay}
does hold, as one sees by integrating with respect to $x''$ for generic fixed $x'$
and using well-known results for oscillatory integrals of the first type
in the sense of \cite{steinbook}.
If $P$ is independent of $x''$ then matters have been reduced to the case where
the ambient dimension is $n$.

As will be clear from its proof,
\eqref{Ltwodecay} continues to hold with uniform constants $C,\eps$ 
if the collection $\{V_j\}$ varies over a compact subset of
the open subset of the relevant Grassmannian manifold consisting
of all such collections in general position, and the polynomials
$P$ are uniformly nondegenerate relative to $\{V_j\}$.

The condition $n<2m$ is necessary for \eqref{Ltwodecay}, with $\|f_j\|_{L^2}$
on the right-hand side rather than $\|f_j\|_{L^\infty}$ as in 
Definition~\ref{defn:powerdecay}.
To see this in the main case where $\{V_j\}$ span $\reals^m$, 
consider without loss of generality the case $n=2m$,
and let each $f_j$ equal the characteristic function of $[-\delta,\delta]$.
If $\delta$ is chosen to be a small fixed constant times $|\lambda|\rp$,
then $|\Lambda_\lambda(f_1,\cdots,f_n)|\ge c\delta^m$.
On the other hand $\prod_j\|f_j\|_{L^2} = c' \delta^{n/2}$.
It follows from the same construction that
if $n>2m$ then the integral defining $\Lambda_\lambda(f_1,\cdots,f_n)$
will in general not even be absolutely convergent for general $f_j\in L^2$.

\begin{theorem} \label{thm:codimone}
Any collection of subspaces $\{V_j\}$ of codimension one 
has the uniform power decay property. 
\end{theorem}

In particular, when $m=2$ then $\kappa=1=m-1$ and hence
nondegeneracy of $P$ implies \eqref{decay}, no matter how large the
collection of subspaces $V_j$ may be.

\begin{theorem} \label{thm:simpledecay} \label{thm:simplequantitative}
Let $m,\kappa$ be arbitrary.
Then any  simply nondegenerate polynomial has the power decay property \eqref{decay}
in every open set. 

More precisely, let $d\in\naturals$ and $c>0$. 
Let $L=\prod_{j=1}^n (w_j\cdot\nabla)$ where each $w_j\in V_j^\perp$ is a  unit vector.
Then there exist $C<\infty$ and $\eps>0$ such that 
for any real-valued polynomial $P$ of degree $\le d$
such that 
$\max_{|x|\le 1}|L(P)(x)|\ge c$,
\begin{equation*}
|\Lambda_\lambda(f_1,\cdots,f_n)|\le C(1+|\lambda|)^{-\eps}\prod_j\|f_j\|_\infty
\end{equation*}
for all functions $f_j\in L^\infty$ and all $\lambda\in\reals$.
\end{theorem}

In each of these theorems, the regularity hypothesis on $\eta$
can be relaxed to H\"older continuity, as discussed above. 

Theorem~\ref{thm:codimone} directly implies Theorem~\ref{thm:simpledecay}.
Indeed, if $w_j$ are as in the latter theorem, define $\tilde V_j = w_j^\perp$,
thus obtaining codimension one subspaces $\tilde V_j\supset V_j$.
Each $f_j\circ\pi_j$ can be rewritten as $\tilde f_j\circ\tilde\pi_j$
where $\tilde\pi_j$ is the orthogonal projection onto $\tilde V_j$
and $\tilde f_j(y)=f_j\circ\pi_j(y)$ for $y\in\tilde V_j$. The hypothesis
$L(P)\ne 0$ guarantees that $P$ is nondegenerate relative to $\{\tilde V_j\}$,
with appropriate uniformity.
On the other hand,
we will show in Proposition~\ref{prop:nondegenequivalence} that for $\kappa=m-1$, 
nondegeneracy is equivalent to simple
nondegeneracy, so that Theorem~\ref{thm:codimone} conversely implies
Theorem~\ref{thm:simpledecay}.

The next two theorems are the best we have been able to do for
general values of the dimension $\kappa$.
\begin{theorem} \label{thm:noslicing}
Let $\kappa\in\{1,2,\cdots,m-1\}$.
Let $\{V_j: 1\le j\le n\}$ be a collection of
$\kappa$-dimensional subspaces of $\reals^m$, lying in general position. If
\[
n<\frac{2m}{\kappa} 
\]
then $\{V_j\}$ has the uniform decay property, in its $\lt$ formulation \eqref{Ltwodecay},
with uniform bounds if $P$ belongs to a family of uniformly nondegenerate polynomials.
\end{theorem}
For $\kappa=m-1$ we have already shown that the conclusion holds with no restriction on $n$.

When the codimension $m-\kappa$ is relatively small but strictly greater than
one then the next result is superior.
\begin{theorem} \label{thm:lowcodimension}
Let $\{V_j: 1\le j\le n\}$ be a finite collection of subspaces of
$\reals^m$, each having dimension $\kappa$ and lying in general position.
If 
\[
n\le \frac{m}{m-\kappa}
\]
then $\{V_j\}$ has the uniform decay property.
\end{theorem}

\subsection{Singular integrals}\label{section:CZO}
\label{subsect:czo}

For any real valued polynomial $P(x,t)$
 of degree  $d$,
consider the operator
\begin{equation}
T(f,g)(x) = \int_{-\infty}^\infty e^{iP(x,t)} f(x-t)g(x+t)\,t\rp\,dt
\end{equation}
where it is initially assumed that $f,g\in C^1_0$,
the class of all continuously differentiable functions having compact supports,
and the integral is taken in the principal-value sense.
One of the purposes
 of this note is to establish the following $L^p$ bounds for these operators.

\begin{theorem}\label{thm1}
For any exponents $p_1,p_2,q\in(0,\infty)$ such that
$q\rp = p_1^{-1}+p_2^{-1}$, $p_1,p_2>1$ and $q>2/3$, and any degree $d\ge 1$,
there exists $C<\infty$
 such that
$\|T(f,g)\|_{q}\le C\|f\|_{p_1}\|g\|_{p_2}$ for all $f,g\in C^1_0$,
uniformly for all real-valued polynomials $P$ of degrees $\le d$.
\end{theorem}

This answers a question raised by Lacey and Thiele.
The cases $d=0,1,2$ were previously known. Indeed,
the case $d=0$ is a celebrated theorem of Lacey and Thiele
\cite{bilinearHilbert}. The case
$d=1$, that is $P(x,t)=a_0x+a_1t$, can be reduced to $d=0$
by replacing $f$ by
$\tilde f(x) = e^{i a_{0}x/2 -ia_1x/2}f(x)$
and $g$ by
$\tilde g(x)=e^{ia_0x/2+ia_1x/2}f(x)$.
The case $d=2$, that is, $P(x,t)=a_0x^2+a_1xt+a_2t^2+b_1x+b_2t$,
 is likewise reducible to $d=0$ by the substitution
$\tilde f(x)=e^{-ia_{1}x^2/4+ia_2x^2/2 }f$,
$\tilde g(x) = e^{ia_1x^2/4+ia_2x^2/2 }g$.
But for $d\ge 3$ no such simple reduction exists.

Theorem~\ref{thm1} is a bilinear analogue of a theorem
of Ricci and Stein \cite{riccistein},
who proved $L^p$ estimates for linear operators
$f\mapsto \int e^{iP(x,t)}f(x-t)K(t)\,dt$,
for arbitrary real-valued polynomials $P$
and Calder\'on-Zygmund kernels $K$.
It is a special case of the following more
general result.
Let $n\geq 1$,
$\Gamma =\{\xi\in \mathbb{R}^{n+1}: \xi_1+\xi_2+\cdots +\xi_{n+1}=0\}$,
and $\Gamma'$ be the orthogonal complement of a subspace of
$
\Gamma$ such that
the dimension of $\Gamma'\cap \Gamma $ is $k$.
Let $K$ be a $k$-dimensional Calder\'on-Zygmund kernel on 
$\Gamma'\cap\Gamma$, that is, $K$ is Lipschitz continuous except at the origin,
$K(r\gamma)\equiv r^{-k}K(\gamma)$ for all $r>0$ and $\gamma\ne 0$,
and $\int_{S^{k-1}}K\,d\sigma=0$ where $\sigma$ denotes
surface measure on the unit sphere $S^{k-1}$ in $\Gamma'\cap\Gamma$.
We define the $n$-linear operator $T$ by
\begin{equation}\label{multi}
T(f_1,f_2, \cdots, f_n)(x)=\int_{\Gamma'\cap\Gamma} f_1(x+\gamma_1)
f_2(x+\gamma_2)\cdots f_n(x+\gamma_n)K(\gamma)d\gamma\,,
\end{equation}
for $x\in\reals^1$,
where $d\gamma $ is Lebesgue measure on $\Gamma'\cap\Gamma$,
and $\gamma_i\in\reals^1$ is the $i$-th coordinate of 
$\gamma\in\reals\times\reals \cdots\times\reals$ 
as an element of $\mathbb{R}^{n+1}$.
The integral is interpreted in the principal-value sense, so that
$T(f_1,\cdots,f_n)(x)$ is well-defined provided that each $f_j\in C^1$
has compact support.
We always  assume that any $k+1$ variables in
$\{x+\gamma_1, \cdots, x+\gamma_n, x\}$ are linearly independent.

To such an operator $T$ and to any real-valued polynomial
$P(x,  \gamma_1,\cdots,\gamma_{n})$ we associate the multilinear
operator
\begin{equation}\label{defofT}
T_P(f_1,f_2,\cdots, f_n)(x)=\int_{\Gamma'\cap\Gamma}
e^{iP(x,\gamma_1,\cdots , \gamma_{n})} f_1(x+\gamma_1)
f_2(x+\gamma_2)\cdots f_n(x+\gamma_n)K(\gamma)d\gamma\,,
\end{equation}
which is again well-defined when each $f_j\in C^1_0$.

Write $p_0 = q'$.
\begin{theorem}  \label{thm2}
Suppose that $n\le 2k$. Then for any
real-valued polynomial $P$, $T_P$ maps
$\otimes_{j=1}^n L^{p_j}$ boundedly to $L^p_0$
whenever $p_0>n^{-1}$, $1<p_j\le\infty$ and
$p_0^{-1} = \sum_j p_j^{-1}$,
provided that
\begin{equation}\label{p4}
\frac{1}{p_{i_1}}+\frac{1}{p_{i_2}}+\cdots +\frac{1}{p_{i_r}} <
\frac{2k+r+1-n}{2}
\end{equation}
for all $0\leq i_1<i_2<\cdots <i_r\leq n$ and all
$1\leq r\leq n+1$.
This conclusion holds uniformly for all polynomials of degrees $\le D$,
for any $D<\infty$. 
\end{theorem}
Under these hypotheses, the nonoscillatory case $P\equiv 0$
was treated in \cite{mtt}.  
Theorems~\ref{thm1} and \ref{thm2} will be proved by combining 
previously known results for nonoscillatory multilinear singular
integral operators with our new results for nonsingular oscillatory
integrals.


\smallskip
\noindent 
{\bf Acknowledgement.}
The second author thanks M.~Lacey for encouraging him to work on this project.

\subsection{Bounds for sublevel sets}
\label{subsect:sublevel}

We are able to prove a weaker consequence of the decay
property in full generality, providing some evidence
that nondegeneracy might always imply it.
For any Lebesgue measurable functions $g_j$ which are finite almost everywhere, 
for any $\eps>0$, and for any bounded subset $B\subset\reals^m$
consider the sublevel sets
\begin{equation}
E_\eps = \big\{
x\in B: 
|P(x)-\sum_j g_j(\pi_j(x))|<\eps
\big\}.
\end{equation}
\begin{theorem}  \label{thm:sublevel}
Suppose that a real-valued polynomial $P$
is nondegenerate with respect to a finite collection
of subspaces $\{V_j\}$ of $\reals^m$.
Then there exists $\delta>0$ so that for any bounded subset $B$ of $\reals^m$
there exists $C<\infty$ such that
for any measurable functions $g_j$ defined on $V_j$,
for any $\eps>0$,
the associated sublevel sets satisfy
\begin{equation} \label{ineq:sublevel}
\big|
E_\eps
\big|
\le C\eps^\delta.
\end{equation}
\end{theorem}

If a real-valued measurable function $P$ has the 
power decay property \eqref{decay}, then the sublevel set bound 
\eqref{ineq:sublevel} holds. Indeed, fix a cutoff function
$0\le h\in C^\infty_0(\reals)$ satisfying $h(t)=1$ whenever
$|t|\le 1$. Fix also $\zeta\in C^\infty_0(\reals^m)$
such that $\zeta\equiv 1$ on $B$.
Then
\begin{align*}
\big|\{x\in B: |P(x)-\sum_j g_j(\pi_j(x))|<\eps\}\big|
&\le \int h[\eps\rp(P-\sum_j g_j\circ\pi_j)(x)]\,\zeta(x)\,dx
\\
&= (2\pi)\rp \eps\int_\reals \widehat{h}(\eps\lambda)
\int e^{i\lambda(P(x)-\sum_j g_j(\pi_j(x))}\zeta(x)\,dx\,d\lambda.
\end{align*}
The inner integral is
$ O(1+|\lambda|)^{-r}$ for some $r\in(0,1)$
by the power decay property applied to $f_j=e^{i\lambda g_j}$.
Inserting this into the double integral gives a majorization
\[
|\{x\in B: |P(x)-\sum_j g_j(\pi_j(x))|<\eps\}|
\le C\int_\reals \eps (1+\eps|\lambda|)^{-2}(1+|\lambda|)^{-r}\,d\lambda
\le C\eps^{-r}.
\]


The remainder of the paper is devoted to proofs of the theorems.
Certain algebraic aspects of nondegeneracy come into play in the analysis
and are developed in \S\ref{section:algebra}; this material is used throughout
the paper. 
The simply nondegenerate case is treated in \S\ref{section:simplynondegenerate}.
The analysis when the dimension $\kappa$ equals $1$ involves an additional ingredient
and is presented in \S\ref{section:oneD}.
Applications to multilinear operators with Calder\'on-Zygmund singularities
and oscillatory factors are carried out in \S\ref{section:czo}.
The bounds for measures of sublevel sets are established in \S\ref{section:sublevel}.
Finally, Theorems~\ref{thm:noslicing} and \ref{thm:lowcodimension} are
proved in the final section.

\section{Algebraic aspects of nondegeneracy}
\label{section:algebra}
This section develops various general properties, largely algebraic in nature,
of nondegeneracy.

As noted above, simple nondegeneracy implies nondegeneracy. In the codimension
one case they are equivalent:
\begin{proposition}  \label{prop:nondegenequivalence}
If each subspace $V_j$ has codimension one, that is if
$\kappa=m-1$, then any nondegenerate polynomial is simply nondegenerate.
\end{proposition}

\begin{proof}
Let a polynomial $P$ and distinct subspaces $V_j$ of codimension $1$ be given.
Let $w_j$ be unit vectors orthogonal to $V_j$,
let $L=\prod_{j=1}^n (w_j\cdot\nabla)$,
and let $L' = \prod_{j=1}^{n-1}(w_j\cdot\nabla)$.
We must show that if $LP\equiv 0$, then $P$ is degenerate.

Proceed by induction on $n$, which is the total number of subspaces $V_j$. 
Then since $L'(w_n\cdot\nabla) P\equiv 0$, 
$(w_n\cdot\nabla) P=\sum_{j=1}^{n-1} q_j\circ\pi_j$ for some polynomials $q_j$.
By a rotation of the coordinate system it may be arranged that $w_n\cdot \nabla = \p_{x_m}$.

$(w_n\cdot\nabla)(f\circ\pi_j) = (D_jf)\circ\pi_j$ for some nonvanishing
constant-coefficient real vector fields $D_j$, for $w_n$ cannot be orthogonal to
$V_j$ for $j<n$ since it is orthogonal to $V_n$, and these codimension one
subspaces are distinct.
Consequently there exist polynomials of the form $p_j\circ\pi_j$ for $1\le j\le n-1$
such that $(w_n\cdot\nabla)(p_j\circ\pi_j) = q_j\circ\pi_j$.

Thus $(w_n\cdot\nabla)(P-\sum_{j<n}p_j\circ\pi_j)\equiv 0$, so the polynomial
$(P-\sum_{j<n}p_j\circ\pi_j)(x)$ is a function of $x_m=\pi_m(x)$ alone,
and hence is of the form $p_n\circ\pi_n$.
This proves that $P$ is degenerate.
\end{proof}


\begin{lemma}  \label{lemma:diffops}
Let $P$ be a homogeneous polynomial of some degree $d$. Then $P$ is nondegenerate
relative to a finite collection of subspaces $\{V_j\}$, of any dimension,
if and only if there exists a constant-coefficient partial differential
operator $\scriptl$, homogeneous of degree $d$, such that $\scriptl(P)\ne 0$
but $\scriptl(p_j\circ\pi_j)\equiv 0$
for every polynomial $p_j:V_j\to\complex$ of degree $d$.
\end{lemma}

\begin{proof}
If such an operator $\scriptl$ exists then $P$ is obviously nondegenerate.
To prove the converse,
fix $d$ and denote by $\scriptp_d$ the
vector space of all homogeneous polynomials of degree $d$.
The pairing
$\langle\scriptl,\,P\rangle = \scriptl(P)$
between 
homogeneous constant-coefficient differential operators
of the same degree $d$ and elements of $\scriptp_d$
is nondegenerate. Thus
the dual space of $\scriptp_d$
may be canonically identified with the vector space of all such differential operators.

If $P$ is nondegenerate then it does not belong to the subspace of all degenerate homogeneous
polynomials, then there exists a linear functional which annihilates that subspace, but
not $P$. This functional may be realized as $Q\mapsto \scriptl(Q)$ for some operator $\scriptl$,
completing the proof in the converse direction.
\end{proof}

Any polynomial may be expressed in a unique way as a sum of homogeneous polynomials
of distinct degrees.
\begin{lemma}
A polynomial $P$ is nondegenerate relative to $\{V_j\}$ if and only if at least one of its
homogeneous summands is nondegenerate. Moreover, a homogeneous polynomial 
is degenerate if and only if it may be expressed as $\sum_j p_j\circ\pi_j$
where each $p_j$ is a homogeneous polynomial of the same degree.
\end{lemma}

\begin{proof}
Both of these assertions follow from the simple fact that a polynomial
$p_j:\reals^\kappa\to\reals$ 
is homogeneous of some degree  if and only if $p_j\circ\pi_j:\reals^k\to\reals$
is homogeneous, of the same degree.
\end{proof}

If $P$ can be decomposed as $\sum_j f_j\circ\pi_j$ for some functions $h_j$
which are not necessarily polynomials,
then the decay bound \eqref{decay} certainly fails, for the same
reason noted earlier (set $f_j = e^{-i\lambda h_j}$). The next lemma says that
it makes no difference whether the $h_j$ are taken to be polynomials,
or arbitrary functions, in the definition of nondegeneracy.

\begin{lemma} \label{lemma:polyordistribution}
A polynomial $P$ is degenerate with respect to a collection of subspaces $\{V_j\}$
if and only if there exist distributions $h_j$ in $\reals^\kappa$ such
that $P=\sum_j h_j\circ\pi_j$ in some open set.
\end{lemma}

\begin{proof}
Let $\{V_j\}$ be given, and let $P$ be any nondegenerate polynomial of some degree $d$.
We need only show that $P$ can't be decomposed locally as $\sum_j h_j\circ\pi_j$,
since the converse is a tautology.
It is no loss of generality
to suppose that the homogeneous summand of $P$ of degree $d$ is nondegenerate.
For otherwise we may express that summand as $\sum_j p_j\circ\pi_j$ where
the $p_j$ are homogeneous polynomials of degree $d$, then replace $P$
by $P-\sum_j p_j\circ\pi_j$ to reduce the degree.

By Lemma~\ref{lemma:diffops}, 
there exists a constant-coefficient linear partial differential operator $\scriptl$,
homogeneous of degree $d$, such that $\scriptl(P)\ne 0$, yet
$\scriptl(p_j\circ\pi_j)=0$  for any polynomial $p_j$ homogeneous
of degree $d$, for any index $j$.

$\scriptl(Q)(0)=0$ for any polynomial $Q:\reals^m\to\reals$ which
is homogeneous of any degree except $d$. By combining this with the hypothesis
we find that 
$\scriptl(p\circ\pi_j)(0)=0$ for any polynomial $p:\reals^\kappa\to\reals$ 
and any index $j$.
Consequently 
$\scriptl(p\circ\pi_j)(x)=0$ for all $x\in\reals^m$.
Since any distribution is a limit of polynomials in the distribution 
topology, 
$\scriptl(h_j\circ\pi_j)$ likewise vanishes identically, in the sense of distributions.
Yet $\scriptl$ does not annihilate $P$,
so $P$ cannot be decomposed as $\sum_j h_j\circ\pi_j$.
\end{proof}

There is a characterization of nondegeneracy in terms of difference operators,
analogous to that involving differential operators in Lemma~\ref{lemma:diffops},
which will be used in the proof of Theorem~\ref{thm:sublevel}.
Let $e_j$ denote the $j$-th coordinate vector in $\reals^m$,
define $\delta_jf(x) = f(x+e_j)-f(x)$,
and define $\Delta_\alpha = \delta_1^{\alpha_1}\circ\cdots\circ
\delta_m^{\alpha_m}$ where $\alpha=(\alpha_1,\cdots,\alpha_m)$
with each $\alpha_j$ an arbitrary nonnegative integer.
Define $|\alpha| = \sum_j \alpha_j$.

If $P:\reals^m\to\reals$ is a homogeneous polynomial of degree $D$
and $|\alpha|=D$, then $\Delta_\alpha(P)$ is a constant.
Thus there is a natural pairing
between such the linear span of such difference operators, and such polynomials. 
This pairing is clearly nondegenerate, 
since whenever $|\alpha|=|\beta|$, $\Delta_\alpha(x^\beta)=0$ if and only 
if $\alpha=\beta$.
Thus it establishes an identification of the span
of all $\Delta_\alpha$ with $|\alpha|=D$
with the dual of the space of all homogeneous real polynomials of degree $D$.

\begin{lemma} \label{lemma:differences}
Let $P:\reals^m\to\reals$ be a polynomial of degree $D$.
Suppose that $P=P_D+R$ where $P_D$ is homogeneous of degree $D$,
$R$ has degree $<D$, and $P_D$ is nondegenerate with respect
to a collection of subspaces $\{V_j\}$.
Then there exist a finite subset $\{y_\beta\}$ of $\reals^m$,
and corresponding scalars $c_\beta\in\reals$,
such that for any $1\le j\le n$ and any continuous function $f_j$,
for all $x$ and all $r>0$,
\begin{equation} \label{zero} 
\sum_\beta c_\beta f_j\circ\pi_j(x+ry_\beta)=0 
\end{equation}
but
\begin{equation}  \label{notzero}
\sum_\beta c_\beta P(x+ry_\beta)\equiv r^D. 
\end{equation}
\end{lemma}

\begin{proof}
By the preceding discussion there exist real scalars $b_\beta$
such that $L=\sum_{|\beta|=D} b_\beta \Delta_\beta$ annihilates
$p\circ\pi_j$ for any homogeneous polynomial $p$ of degree $D$,
for all $j$,
but $L(P_D)\ne 0$.
$L$ of course annihilates $p\circ\pi_j$ whenever $p$
has lower degree, and $L(p\circ\pi_j)(0)=0$ whenever $p$ is
homogeneous of some higher degree.
It follows that $L(\sum_j p_j\circ\pi_j)(x)=0$ for
all $x\in\reals^m$, for all polynomials $p_j$,
and from this the same follows for all continuous functions
$p_j$; it even holds in the sense of distributions for arbitrary
distributions $p_j$.
Since $R$ has degree $<D$, $L$ also annihilates $R$,
so $L(P)$ is a nonzero constant.
By multiplying by a suitable scalar, we may arrange for this constant
to be $1$.

$L$ takes the form $L(Q)(x) = \sum_\beta c_\beta Q(x+y_\beta)$.
Since $P_D$ is homogeneous and any translate and dilate of $R$
has degree $<D$, \eqref{notzero} follows.
\eqref{zero} for general $r$ follows from the case $r=1$
by scaling.
\end{proof}

We say that a polynomial is a {\em monomial}
if it can be expressed as a product of linear factors.
By a {\em differential monomial} we mean a linear partial differential
operator which can be expressed as a product of finitely many real vector fields,
each with constant coefficients.

\begin{lemma} \label{lemma:lowcodimension}
Let $\{V_j: 1\le j\le n\}$ be a finite collection of subspaces of
$\reals^m$, each having dimension $\kappa$ and lying in general position.
If $(m-\kappa)n\le m$ then any polynomial which is nondegenerate
relative to $\{V_j\}$ is simply nondegenerate relative to
$\{V_j\}$.
\end{lemma}

\begin{proof}
Consider the variety $\scriptv = \cup_j\scriptv_j$
where $\scriptv_j\subset\reals^m$ is the $m-\kappa$-dimensional subspace of $\reals^m$
equal to the range of $\pi_j^*$, where $\pi_j^*$ is the adjoint of
the orthogonal projection $\pi:\reals^m\to V_j$.
Let $\scripti$ be the ideal in $\complex[x]$
consisting of all real polynomials which vanish on $\scriptv$.
Then $\scripti$ is of course finite-dimensional as a $\complex[x]$-module.
According to the next lemma,
the ideal $\scripti$ is generated by some finite set of monomials.

Now consider the ring $\scriptl$
of all constant-coefficient linear partial differential
operators in $\reals^m$, and in it the ideal
$\scriptl_{\{V_j\}}$
consisting of all elements which annihilate all polynomials
$p_j\circ\pi_j$, for all $1\le j\le n$.
The map $L\mapsto\sigma(L)$ from $L$ to its symbol is an isomorphism
of $\scriptl$ with $\complex[x]$, which maps
$\scriptl_{\{V_j\}}$ to $\scripti$.
Thus any element $L\in \scriptl_{\{V_j\}}$
may expressed as a finite sum $\sum_\alpha \ell_\alpha\circ
\scriptm_\alpha$ where each $\scriptm_\alpha$ is a differential monomial
which likewise annihilates all $p_j\circ\pi_j$,
and the coefficients $\ell_\alpha$ belong to $\scriptl$.

If a polynomial $P$ is nondegenerate relative to $\{V_j\}$,
then there exists $L\in \scriptl_{\{V_j\}}$ such that $L(P)$
does not vanish identically. Therefore there exists a differential
monomial $\scriptm_\alpha$ such that $\scriptm_\alpha(P)$ does not vanish identically,
so that $P$ is simply nondegenerate.
\end{proof}

\begin{lemma}
Let $\scriptv = \cup_{j=1}^n \scriptv_j$ where each $\scriptv_j$
is a linear subspace of $\complex^m$ of dimension $\kappa$.
Suppose that $(m-\kappa)n\le m$, and that the subspaces
$\scriptv_j$ lie in general position.
Then the ideal $\scripti\subset\complex[x]$ of all polynomials
vanishing identically on $\scriptv$ is generated by monomials.
\end{lemma}

\begin{proof}
Fix linear functions $y_{j,k}(x)$, for $1\le j\le n$
and $1\le k\le m-\kappa$, such that
$\scriptv_j=\{x: y_{j,k}(x)=0 \text{ for all } 1\le k\le m-\kappa\}$.
The general position hypothesis means simply that these
are linearly independent over $\complex$.
Thus after possibly adding some additional coordinates if $(m-\kappa)n<m$,
we may regard these functions as coordinates on $\complex^m$.
Any monomial of the form $\prod_{j=1}^n y_{j,k(j)}$ belongs to $\scripti$. 

Let $P\in\scripti$. We wish to show that $P$  belongs to the ideal
generated by all products 
$\prod_{j=1}^n y_{j,k(j)}$,
where the function $j\mapsto k(j)$ runs over all $n^{m-\kappa}$ possibilities.

Since $P\equiv 0$ on $\scriptv_1$,
it may be expressed as $P = \sum_{k=1}^{m-\kappa} y_{1,k}r_{1,k}$
for certain polynomials $r_{1,k}$, which have the additional property
that $r_{1,k}$ is independent of $y_{1,i}$ for all $i<k$.
Consider the restriction of $P$ to the subspace where $y_{1,k}=0$
for all $k<m-\kappa$. Since $P$ vanishes on the intersection
of this subspace with $\scriptv_2$, so must $y_{1,m-\kappa}r_{1,m-\kappa}$,
and hence, by the general position hypothesis, so must $r_{1,m-\kappa}$.
Therefore $r_{1,m-\kappa}$ must belong to the ideal generated
by $\{y_{2,k}\}$.

Next consider the restriction of $P$ to the subspace where
$y_{1,k}=0$ for all $k<m-\kappa-1$.
$r_{1,m-\kappa}\equiv 0$ on $\scriptv_2$, so by
repeating the reasoning of the preceding paragraph we may conclude that
$r_{1,m-\kappa-1}$ belongs to the ideal generated by $\{y_{2,k}\}$.
By descending induction on $k$ we eventually find that every $r_{1,k}$
does also.

Thus we may write $P = \sum_{k_1=1}^{m-\kappa}\sum_{k_2=1}^{m-\kappa}
y_{1,k_1}y_{2,k_2}r_{1,k_1,2,k_2}$
for certain polynomials $r_{1,k_1,2,k_2}$.
By rearranging these we may also ensure that $r_{1,k_1,2,k_2}$
is independent of $y_{1,j}$ for all $j<k_1$, and of $y_{2,j}$
for all $j<k_2$.

The same reasoning may now be applied to $\scriptv_3$,
and by induction to $\scriptv_j$ for all $j\le n$.
\end{proof}

\begin{example}
In $\reals^4$ let 
$\scriptv_1=\{x: x_1=x_2=0\}$,
$\scriptv_2=\{x: x_3=x_4=0\}$,
$\scriptv_3=\{x: x_1-x_3=x_2-x_4=0\}$.
Thus $m=4$, $\kappa=2$, and $n=3>m/(m-\kappa)$.
These three subspaces lie in general position.
The polynomial $p(x) = x_1x_4-x_2x_3$ belongs to the associated ideal $\scripti$.
($p$ vanishes on other two-dimensional subspaces as well, but these together
with the three subspaces listed do not lie in general position.)
On the other hand, any monomial $\scriptm$ which belongs to $\scripti$
must have degree at least three.
Indeed, if $\scriptm$ has degree two then
in order to vanish on $\scriptv_1\cup\scriptv_2$, it must
take the form $\scriptm(x)=(ax_1+bx_2)(cx_3+dx_4)$. 
Then $\scriptm(x_1,x_2,x_1,x_2)$ plainly cannot vanish identically on $\scriptv_3$. 
\end{example}

\section{The simply nondegenerate case}
\label{section:simplynondegenerate}

In proving \eqref{decay}, 
each function $f_j$ may be assumed to be supported in
a fixed compact set, the image under $\pi_j$ of the support of the cutoff
function $\eta$. We will assume this throughout the proofs of
Theorems \ref{thm:slices} and \ref{thm:simpledecay}. 

\begin{proof}[Proof of Theorem~\ref{thm:simplequantitative}]
The proof proceeds by induction
on the degree $n$ of multilinearity of $\Lambda_\lambda$, and  
is an adaptation of an argument in \cite{ccw}.
Adopt coordinates in which $V_n=\{(y,z)\in\reals^{m-1}\times\reals^{1}: z=0\}$
and hence $\pi_n(y,z) =y$.
Write
\[
\Lambda_\lambda(f_1,\cdots,f_n) = \int f_n(y) \big(\int e^{i\lambda P(y,z)}
\prod_{j=1}^{n-1} f_j(\pi_j(y,z))\,\eta(y,z)\,dz \big)\,dy.
\]
This equals $\langle T_\lambda(f_1,\cdots,f_{n-1}),\overline{f_n}\rangle$ 
for a certain linear operator $T_\lambda$, 
whence 
\[
|\Lambda_\lambda(f_1,\cdots,f_n)|\le \|f_n\|_2
\|T_\lambda(f_1,\cdots,f_{n-1})\|_2\le C \|f_n\|_\infty \|T_\lambda(f_1,\cdots,f_{n-1})\|_2.
\]
Thus it suffices to bound $T_\lambda$ as an operator from $L^\infty\times\cdots\times L^\infty$
to $L^2$. 

Now $\int |T_\lambda(f_1,\cdots,f_{n-1})(y)|^2\,dy$ equals
\[
\int_{\reals^{m-1}}
\iint_{\reals^2} e^{i\lambda(P(y,z)-P(y,z'))}
\prod_{j<n} f_j(\pi_j(y,z))\overline{f_j}(\pi_j(y,z'))
\eta(y,z)\overline{\eta}(y,z')\,dz\,dz'\,dy.
\]
Define $Q_\zeta(y,z) = P(y,z)-P(y,z+\zeta)$, and
$\tilde\eta_\zeta(y,z) = \eta(y,z)\overline{\eta}(y,z+\zeta)$.
Likewise define $F_j^\zeta: V_j\to\complex$ so that
$F_j^\zeta\circ\pi_j(y,z) = f_j(\pi_j(y,z))\overline{f_j}(\pi_j(y,z+\zeta))$;
the right-hand side is a function of $(\pi_j(y,z),\zeta)$ alone because
of the linearity of $\pi_j$.
Of course $\|F_j^\zeta\|_\infty\le \|f_j\|_\infty^2$.
With these definitions and the substitutions $x=(y,z)$, $z' = z+\zeta$,
\[
\int |T_\lambda(f_1,\cdots,f_{n-1})(y)|^2\,dy
= \int \Lambda_\lambda^\zeta (F_1^\zeta,\cdots,F_{n-1}^\zeta)\,d\zeta
\]
where
\[
\Lambda_\lambda^\zeta (F_1^\zeta,\cdots,F_{n-1}^\zeta)
= \int_{\reals^m} e^{i\lambda Q_\zeta(x)} \prod_{j=1}^{n-1}F_j^\zeta(\pi_j(x))
\tilde\eta_\zeta(x)\,dx.
\]
The outer integral may be taken over a bounded subset of $\reals^{m-\kappa}$,
since $\tilde\eta_\zeta\equiv 0$ when $|\zeta|$ is sufficiently large.
Henceforth we assume $\zeta$ to be restricted to such a bounded set.

For each $\zeta$, consider the polynomial phase $Q_\zeta$. 
We may assume that $|\lambda|\ge 1$, since 
\eqref{decay} holds trivially otherwise.
It is given as a hypothesis that there exists a differential operator
of the form $L=\prod_{j=1}^n (w_j\cdot\nabla)$,
with each $w_j\in V_j^\perp$, such that 
$\sup_{|x|\le 1} |L(P)(x)|\ge c>0$.
Let $L' = \prod_{j<n} (w_j\cdot\nabla)$. 
For any $\rho\in(0,1)$ define
\[
E_\rho
=\{\zeta: \max_{|x|\le 1} |L'Q_\zeta(x)|\le\rho\}.
\]

For any $\zeta\notin E_\rho$,
write $\lambda Q_\zeta = (\lambda \rho) \tilde Q$
where $\tilde Q = \rho\rp Q_\zeta$.
Then by the induction hypothesis, there exist $C,\eps'\in\reals^+$ such that 
for all $\zeta\notin E_\rho$,
\begin{equation*}
|\Lambda_\lambda^\zeta (F_1^\zeta,\cdots,F_{n-1}^\zeta)|
\le C(1+|\lambda|\rho)^{-\eps'}\prod_{j<n}\|F_j^\zeta\|_\infty
\le C(1+|\lambda|\rho)^{-\eps'}\prod_{j<n}\|f_j\|_\infty^2.
\end{equation*}
The same bound holds for the integral over all $\zeta\notin E_\rho$,
since $\zeta$ is confined to a bounded set.

For $\zeta\in E_\rho$ there is the trivial estimate
\[
|\Lambda_\lambda^\zeta (F_1^\zeta,\cdots,F_{n-1}^\zeta)|\le C\prod_{j<n}\|F_j^\zeta\|_\infty,
\]
so
\[
\int_{E_\rho} \Big|\Lambda_\lambda^\zeta (F_1^\zeta,\cdots,F_{n-1}^\zeta) \Big|\,d\zeta
\le C|E_\rho| \prod_{j<n} \|f_j\|_{L^\infty}^2.
\]

Now
\begin{equation} \label{sublevelvolume}
|E_\rho|\le C\rho^\delta
\ \text{ for some $\delta>0$ and $C<\infty$.} 
\end{equation}
Indeed, by hypothesis
$\sup_{(x,\zeta)} |\partial_\zeta(L'Q_\zeta(x))|\ge c>0$, and as is well known,
this implies that if $(x,\zeta)$ is restricted to any fixed bounded set, then 
\begin{equation}
|\{(x,\zeta): |L'Q_\zeta(x)|\le r\}|\le Cr^a
\end{equation}
for some $a>0$ where $C,a$ depend only on $c$ and on an upper bound for the degree
of $Q$ as a polynomial in $(x,\zeta)$.
In particular, so long as $x,\zeta$ are restricted to lie in any fixed bounded set,
\[
|\{\zeta: \max_x |L'Q_\zeta(x)|\le r\}|\le Cr^a.
\]

Thus in all 
\[
\int\big|\Lambda_\lambda^\zeta(F_1^\zeta,\cdots,F_{n-1}^\zeta) \big|\,d\zeta
\le C[(|\lambda|\rho)^{-\eps'} + \rho^\delta] \prod_{j<n}\|f_j\|_\infty^2.
\]
Choosing $\rho= |\lambda|^{-c}$ for any fixed $c\in(0,1)$ yields the desired bound.
\end{proof}

\section{The power decay property for $\kappa=1$}
\label{section:oneD}


This section is devoted to the proof of Theorem~\ref{thm:slices}. The proof
turns on a concept related to a notion of uniformity
employed by Gowers \cite{gowers}.
Let $d\ge 1$, and fix a bounded ball $B\subset\reals^m$.
Let $\tau>0$ be a small quantity to be chosen below, let $|\lambda|\ge 1$,
and consider any function $f\in \lt(\reals^m)$ supported in $B$.
\begin{definition} \label{defn:uniformity}
$f$ is $\lambda$-nonuniform if
there exist a polynomial $q$ of degree $\le d$ and a scalar $c$ such that
\begin{equation} \label{nonuniform}
\|f-ce^{iq}\|_{L^2(B)} \le (1-|\lambda|^{-\tau})\|f\|_{L^2}.
\end{equation}
Otherwise $f$ is said to be $\lambda$-uniform.
\end{definition}

This notion depends on the parameters $d,\tau$.
So long as they remain fixed
there exists $C<\infty$, depending also on $B$,
such that any $\lambda$-uniform function $f\in \lt(B)$
satisfies favorable bounds for generalized Fourier coefficients:
\begin{equation} \label{gowersuniformity}
\big|\int f(t)e^{-iq(t)}\,dt\big|\le C|\lambda|^{-\tau/2}\|f\|_{L^2(B)} 
\end{equation}
unifomly for all real-valued polynomials $q$ of degree $\le d$. Indeed, 
$f$ could otherwise be decomposed  in $\lt(B)$
into its projection onto $e^{iq}$ plus an orthogonal vector, implying \eqref{nonuniform}.

The proof of Theorem~\ref{thm:slices} will proceed by induction on $n$, 
the number of subspaces $\{V_j\}$.
The inductive step enters in the following way: If $f_1 = e^{ip}$
for some polynomial $p$,
then $\Lambda_\lambda(e^{ip},f_2,\cdots,f_n)
=\int e^{i\tilde P(x)}\prod_{j=2}^n f_j(\pi_j(x))\,\eta(x)\,dx$
where $\tilde P = P+p\circ\pi_1$. With $p$ held fixed,
this may be regarded as an $(n-1)$-multilinear operator acting
on $(f_2,\cdots,f_n)$, with the new phase $\tilde P$ and
the smaller collection of subspaces $\{V_j: 2\le j\le n\}$. 
$\tilde P$ is nondegenerate relative to this subcollection; moreover,
this nondegeneracy is uniform as $p$ varies over all polynomials of
uniformly bounded degrees.
Thus it is a consequence of the inductive hypothesis that
\begin{equation} \label{inductive}
|\Lambda_\lambda(e^{ip},f_2,\cdots,f_n)|\le C|\lambda|^{-\eps}
\end{equation}
provided that $\|f_j\|_2\le 1$ for all $2\le j\le n$,
uniformly for all polynomials $p$ of uniformly bounded degrees.

\begin{proof}[Proof of Theorem~\ref{thm:slices}]

Given the cutoff function $\eta$, there exist intervals $B_j$ of finite lengths
in $\reals^1$ such that $\Lambda_\lambda(f_1,\cdots,f_n)$
depends only on the restriction of each $f_j$ to $B_j$,
so henceforth we assume $f_j$ to be supported in $B_j$.
Thus $\|f_j\|_{\lt}$ equals $\|f_j\|_{\lt(B_j)}$.

We may assume henceforth\footnote{Our proof can be
adapted to the simpler case $n\le m$ as well.}  
that $n$ is strictly larger than $m$,
since the theorem is already known in a more precise form 
\cite{phongsteinsturm},\cite{riccistein},\cite{steinbook}  for the case $n=m$.
Let $e_1$ be a unit vector orthogonal to the span of $\{V_j: 2\le j\le m\}$.
Since the subspaces $\{V_j: 2\le j\le m\}$ span a space of codimension one,
$e_1$ is uniquely determined modulo multiplication by $-1$, and it cannot
be orthogonal to $V_1$ because of the general position hypothesis.

Likewise choose some unit vector
$e_2$ orthogonal to $\Span(\{V_j: j>m\})$, and not orthogonal to $V_1$.
$\cup_{j>m}V_j$ spans a space
of dimension $n-m <m$ by the assumption that $n<2m$ and the general position
hypothesis, and the only way that all unit vectors in its orthocomplement
could be forced to be orthogonal to $V_1$ is if $V_1$ were to be contained
in $\Span(\{V_j: j>m\})$. But this cannot happen, again by general position
and the restriction $n-m<m$. Thus there exists at least one vector $e_2$ with the
required properties.

$e_2$ cannot be orthogonal to $\Span(\{V_j: 2\le j\le m\})$, 
since $\cup_{j=2}^n V_j$ spans $\reals^m$ by the general position assumption 
($n>m$).
$e_1$ cannot be orthogonal to $\Span(\{V_j: j>m\})$, since
then it would be orthogonal to $\cup_{j\ge 2}V_j$, which we have just seen
to be impossible. Thus
$e_2$ is automatically linearly independent of $e_1$.

Given $\{V_j\}$, $P$, a cutoff function $\eta$, and $\lambda$,
define $A(\lambda)$ to be the best constant (that is, the infimum of all
admissible constants) in the inequality
\begin{equation}
|\Lambda_\lambda(f_1,\cdots,f_n)|\le A(\lambda)\prod_j\|f_j\|_{L^2}.
\end{equation}
As noted before the statement of Theorem~\ref{thm:slices}, the hypotheses of general
position and $n\le 2m$ ensure that $A_\lambda$ is finite for all $\lambda$.
We will assume that $\|f_j\|_{L^2}\le 1$ for
all $j$, and that $|\lambda|$ exceeds some sufficiently large constant. 
To prove the theorem it suffices to obtain 
a bound of $C|\lambda|^{-\eps}$ under these additional hypotheses.

The analysis of $\Lambda_\lambda(f_1,\cdots,f_n)$ is divided into two 
cases, depending on whether or not $f_1$ is $\lambda$-uniform. 
If it is not, let $f=f_1,c,q$ satisfy \eqref{nonuniform}. Then 
\[
|\Lambda_\lambda(f_1-ce^{iq},f_2,\cdots,f_n)|\le A(\lambda)(1-|\lambda|^{-\tau}).
\]
Moreover $|c|$ is majorized by an absolute constant since $\|f_1\|_{L^2}=1$
and $e^{iq}$ is unimodular,
so by the inductive hypothesis \eqref{inductive} 
\[
|\Lambda_\lambda(ce^{iq},f_2,\cdots,f_n)|\le C|\Lambda_\lambda(e^{iq},f_2,\cdots,f_n)|
\le C|\lambda|^{-\sigma}
\]
for certain $C,\sigma\in(0,\infty)$; this bound holds uniformly provided that
$P$ varies over a set of uniformly nondegenerate polynomials.
Combining the two terms yields
\begin{equation} \label{nonuniformconclusion}
|\Lambda_\lambda(f_1,\cdots,f_n)|\le
 A(\lambda)(1-|\lambda|^{-\tau}) + |\lambda|^{-\sigma}
\end{equation}
provided that $\|f_j\|_{\lt}\le 1$ for all $j$, 
and that $f_1$ is $\lambda$-nonuniform.
 
Suppose next that $f_1$ is $\lambda$-uniform.
Adopt coordinates $(t,y)\in\reals^2\times\reals^{m-2}$
where $\reals^2$ is the span of $e_1,e_2$ and $\reals^{m-2}$ is its
orthocomplement, by writing $x = t_1e_1+t_2e_2+y$.
Let $P^y(t) = P(t,y)$.
Define 
\[F_1^y(t_2) = \prod_{j=2}^m f_j\circ\pi_j(t,y);\]
the right-hand side is independent of $t_1$ since $e_1$ is orthogonal
to $V_j$ for all $2\le j\le m$.
Likewise define
\[F_2^y(t_1) = \prod_{j=m+1}^n f_j\circ\pi_j(t,y),\]
which is independent of $t_2$.
In the same way, since neither $e_1$ nor $e_2$ is orthogonal to $V_1$,
$f_1\circ\pi_1(t,y)$ depends for each $y$ only on a certain projection
$\pi(t)$ of $t\in\reals^2$ onto $\reals^1$, so we may  write
\[G^y(\pi(t))= f_1\circ\pi_1(t,y).\] 

The general position hypothesis, the hypothesis $n\le 2m$,  and Fubini's theorem 
together imply that
\begin{equation*} 
\begin{aligned}
&\int \|F_1^y\|_{L^2}^2\|G^y\|_{\lt}^2\,dy = C\prod_{j=1}^m \|f_j\|_{L^2}^2
\\
&\int \|F_2^y\|_{L^2}^2\,dy = C\|f_1\|_{L^2}^2\prod_{j=m+1}^n \|f_j\|_{L^2}^2,
\end{aligned}
\end{equation*}
and consequently
\begin{equation} \label{sliceL2bound} 
\int \|F^y_1\|_{L^2}\|F^y_2\|_{L^2}\|G^y\|_{L^2}\,dy
\le C\prod_{j=1}^n \|f_j\|_{L^2}
\end{equation}
by Cauchy-Schwarz.

For each $y\in\reals^{m-2}$ consider
\[
\Lambda_\lambda^y = \int_{\reals^2} e^{i\lambda P^y(t)}F^y_1(t_2)F^y_2(t_1)G^y(\pi(t))\eta(t,y)\,dt.
\]
Define the set $\scriptb\subset\reals^{m-2}$ of bad parameters to 
be the set of all $y$ for which $P^y$ may be decomposed as
\[P^y(t)=Q_1(t_1)+Q_2(t_2)+Q_3(\pi(t))+R(t)\]
where $Q_j$ are real-valued polynomials of degrees $\le d$ on $\reals^1$ for $j=1,2$,
$t_j$ is a certain linear function of $t$,
$R:\reals^2\to\reals$ is likewise a real-valued polynomial of degree $\le d$,
and 
\[\|R\|\le|\lambda|^{-\rho};\]
that is, $P^y$ has small norm
in the quotient space of polynomials modulo degenerate polynomials,
relative to the three projections $t\mapsto t_1,t_2,\pi(t)$. 
Here $\|\cdot\|$ denotes any fixed norm
on the vector space of all polynomials of degree $\le d$, 
and $\rho\in(0,1)$ is another parameter to be specified.
If $y\notin\scriptb$ then $|\lambda|^\rho P^y$ is at least a fixed positive
distance from the span of all polynomials $Q_1(t_1)+Q_2(t_2)+Q_3(\pi(t))$,
so we may apply Theorem~\ref{thm:simplequantitative} 
with $n=3$, $m=2$, and the phase $|\lambda|^\rho P^y$, to obtain 
\begin{equation} \label{goodparametercase}
|\Lambda_\lambda^y|\le C(|\lambda|^{1-\rho} )^{-\tilde\rho}\|F_1^y\|_{L^2}
\|F^y_2\|_{L^2}\|G^y\|_{L^2}
\ \text{for some $C,\tilde\rho\in\reals^+$.}
\end{equation}
This together with \eqref{sliceL2bound} implies that
\[
\int_{y\notin\scriptb} |\Lambda_\lambda^y|\,dy \le C|\lambda|^{-(1-\rho)\tilde\rho},
\]
as desired.

Bad parameters cannot be handled in this way. Nor is it true that
the set of bad parameters has measure $O(|\lambda|^{-\delta})$ for some $\delta>0$;
it could happen that every parameter is bad. Nonetheless we claim that 
if the exponent $\rho$ appearing in the definition of bad parameters
is chosen to be sufficiently small, then there exists
$\eps>0$ such that uniformly for all $y\in\scriptb$,
\begin{equation} \label{keyclaim}
|\Lambda_\lambda^y|\le C|\lambda|^{-\eps}\|F^y_1\|_{L^2}\|F^y_2\|_{L^2},
\end{equation}
which by \eqref{sliceL2bound} again implies the desired bound.
To verify this fix $y\in\scriptb$, and let $P^y(t) = Q_1(t_1)+Q_2(t_2)+Q_3(\pi(t))+R(t)$
as above. Write $\pi(t) = c_1t_1+c_2t_2$ where 
$c_1,c_2$ are both nonzero. Introduce the functions
\begin{align*}
\tilde F_1(t_2) &= F_1^y(t_2)e^{i\lambda Q_2(t_2)},
\\
\tilde F_2(t_1) &= F_2^y(t_1)e^{i\lambda Q_1(t_1)},
\\
\tilde G(s) &= G^y(s) e^{i\lambda Q_3(s)},
\\
\zeta(t) &= \eta(t,y)e^{i\lambda R(t)}.
\end{align*}
Since $\lambda R$ is a polynomial of bounded degree which is $O(|\lambda|^{1-\rho})$
on the support of $\zeta$, the same holds for all its derivatives and therefore
\[
|\widehat{\zeta}(\xi)|\le C_N|\lambda|^{1-\rho}(1+|\xi|)^{-N} \ \text{for all $\xi$ and all $N$.}
\]

By Fourier inversion, 
\[
\Lambda_\lambda^y = c\int_{\reals^3} 
\widehat{\tilde G}(\xi_0)\widehat{\tilde F_1}(-c_2\xi_0-\xi_2)
\widehat{\tilde F_2}(-c_1\xi_0-\xi_1)\widehat{\zeta}(\xi_1,\xi_2)\,d\xi_0d\xi_1d\xi_2.
\]
From this representation, the bound for $\widehat{\zeta}$, and Plancherel's theorem
it follows that
\[
|\Lambda_\lambda^y| \le C|\lambda|^{1-\rho}
\|\widehat{\tilde G}\|_{L^\infty} \|\tilde F_1\|_{L^2}\|\tilde F_2\|_{L^2}.
\]
Here $\|\tilde F_i\|_{L^2} = \|F^y_i\|_{L^2}$.
$\tilde G$ is obtained from $f_1$ by composing with a fixed linear mapping,
translating by an arbitrary amount, and multiplying with $e^{iq}$ for some
real-valued polynomial $q$ of degree $ \le d$. Therefore \eqref{gowersuniformity}
implies that 
\[\|\widehat{\tilde G}\|_{L^\infty}\le C|\lambda|^{-\tau}\]
for some $C,\tau\in\reals^+$ independent of $y$. 
Thus if $f_1$ is $\lambda$-uniform then 
\[|\Lambda_\lambda^y| \le C|\lambda|^{1-\rho-\tau} \|F^y_1\|_{L^2} \|F^y_2\|_{L^2}
\le C |\lambda|^{1-\rho-\tau}.\]
The parameter $\rho\in(0,1)$ is at our disposal; we choose $\rho<1$ sufficiently
close to $1$ such that $1-\rho-\tau=-\eps$ is strictly negative. 
This completes the analysis of the case where $f_1$ is $\lambda$-uniform.

Combined with \eqref{nonuniformconclusion} and \eqref{goodparametercase},
this proves that 
\begin{equation}
A(\lambda)\le C
\max\big(|\lambda|^{-\eps}, |\lambda|^{-(1-\rho)\tilde\rho},
 A(\lambda)(1-|\lambda|^{-\tau}) + |\lambda|^{-\sigma}) \big).
\end{equation}
Therefore $A(\lambda)$ is majorized by a constant times some negative
power of $|\lambda|$, as was to be proved.
\end{proof}

\section{Analysis of oscillatory singular integral operators}
\label{section:czo}
In this section we combine the results on nonsingular oscillatory integral operators proved above
with the currently existing theory of multilinear Calder\'on-Zygmund singular
integral operators to prove Theorems~\ref{thm1} and \ref{thm2}.
Let $d,k,n\ge 1$, and
let $K$ be a Calder\'on-Zygmund distribution in $\reals^{k}$. Thus $K$ is a tempered
distribution which agrees with a Lipschitz function on $\reals^k\setminus\{0\}$,
and $K(ry)\equiv r^{-k}y$ in the sense of distributions.
For $1\le j\le n$ let $\ell_j:\reals^k\to\reals^d$ be a linear mapping, and
suppose that the intersection of the nullspaces of all the $\ell_j$ is $\{0\}$.
Consider the multilinear operators
\begin{equation}
T(f_1,\cdots,f_n)(x) = \big\langle K(y),\ \prod_{j=1}^n f_j(x+\ell_j(y))\big\rangle
\end{equation}
where $f_j\in C^\infty_0$, and the pairing is that of the distribution $K$ 
with the test function $y\mapsto \prod_{j=1}^n f_j(x+\ell_j(y))$.
To any real-valued polynomial $P(x,y)$ defined on $\reals^{d+k}$ is associated
the operator
\begin{equation}
T_P(f_1,\cdots,f_n)(x) = \big\langle e^{iP(x,y)}K(y),\ \prod_{j=1}^n f_j(x+\ell_j(y))\big\rangle.
\end{equation}

It is not known under what conditions operators $T$
act boundedly on Lebesgue spaces, but we can assert 
a conditional result, generalizing Theorem~\ref{thm2}.
Set $\ell_0(y)\equiv 0$, and in $\reals^{d+k}$ let $V_j$ be the orthocomplement
of the nullspace of the mapping $(x,y)\mapsto x+\ell_j(y)$.
Thus $x+\ell_j(y)$ may alternatively be written as $\pi_j(x,y)$.

\begin{theorem}  \label{thm:conditional}
Let $d,n,n,K,T,\{\ell_j\}$ be as above.
Suppose that $T$ maps $\otimes_{j=1}^n L^{p_j}$
boundedly to $L^q$ for some exponents $1<p_j<\infty$
and $0<q<\infty $ which satisfy $q^{-1} = \sum_j p_j^{-1}$.
Suppose furthermore that the associated subspaces $\{V_j: 0\le j\le n\}$
have the uniform decay property. 
Then for any real-valued polynomial $P$,
$T_P$ maps $\otimes_{j=1}^n L^{p_j}$
boundedly to $L^q$. Moreover  for any $D$ there exists $C_D<\infty$
such that
$\|T_P(f_1,\cdots,f_n)\|_q\le CC_D\prod_j\|f_j\|_{p_j}$
uniformly for all real-valued polynomials $P$ of degrees $\le D$,
where $C$ is the operator norm of $T$ from $\otimes_j L^{p_j}$ to $L^q$.
\end{theorem}

If $\{V_j\}$ merely has the uniform decay property for all polynomials of
degrees $\le D$, then the conclusion of the theorem holds for
polynomials of degrees $\le D$.

A sufficient condition for $T$ to be bounded was established
in \cite{mtt}. Combining it with Theorem~\ref{thm:conditional} gives
Theorem~\ref{thm2}, and its special case Theorem~\ref{thm1}.

The proof of Theorem~\ref{thm:conditional} is little different
from its linear prototype \cite{riccistein}. One step in the proof 
requires passing to associated truncated operators. This can be done quite generally
by means of the following lemma, versions of which have appeared elsewhere
in the literature.

\begin{lemma}\label{lemma:truncation}
Let $m\ge\kappa\ge 1$ and $n\ge 0$ be integers.
For $0\le j\le n$ let $\ell_j:\reals^m\to\reals^\kappa$ be linear mappings, 
the intersection of all of whose nullspaces is $\{0\}$.
Let $p_j\in[1,\infty]$. 

Let $\eta\in C^{m+1}_0(\reals^m)$ have compact support.
There exists $A<\infty$, depending only on the $C^{m+1}$ norm and support of $\eta$,
with the following property.
Let $\scriptk$ be any tempered distribution in $\reals^m$,
and suppose there exists $C<\infty$ such that 
\begin{equation} \label{lemmahypothesis}
\langle \scriptk,\,\prod_{j=0}^n f_j\circ\ell_j\rangle|\le C\prod_j\|f_j\|_{p_j}
\ \ \text{for all $f_j\in C^\infty_0(\reals^\kappa)$.}
\end{equation}
Then 
\begin{equation} \label{lemmaconclusion}
|\langle \eta \scriptk,\,\prod_{j=0}^n f_j\circ\ell_j\rangle|\le AC\prod_j\|f_j\|_{p_j}
\ \ \text{for all $f_j\in C^\infty_0(\reals^\kappa)$.}
\end{equation}
\end{lemma}
\noindent
Here $\langle \scriptk,\,\cdot\rangle$ denotes the pairing of distribution with test function;
the hypothesis that the intersection of the nullspaces is trivial guarantees
that $\prod_{j=0}^n f_j\circ\ell_j$ is compactly supported.

\begin{remark}
The corresponding conclusion holds for multilinear {\em operators}
mapping $\otimes_{j=1}^n L^{p_j}$ to $L^q$, for any $q\in(0,\infty]$,
with the sole modification that when $q<1$, the hypothesis $\eta\in C_0^{m+1}$
should be changed to $\eta\in C_0^s$ for some suitably large $s=s(q,m)$. 
This follows from the argument below, together with
the quasi-triangle inequality $\|f+g\|_q^q\le \|f\|_q^q+\|g\|_q^q$.
\end{remark}

\begin{proof}
Write $\scriptt_\scriptk(f_0,\cdots,f_n) = \langle \scriptk,\,\prod_j f_j\circ\ell_j\rangle$.
Fix $\zeta\in C^\infty_0(\reals^\kappa)$ such that
$\zeta\circ\ell_j(y)\equiv 1$ for all $y$ in the support of $\eta$.
Thus
\[
\scriptt_{\eta \scriptk}(f_0,\cdots,f_n)
=
\langle \scriptk\eta,\,\prod_{j=0}^n (f_j\cdot \zeta)\circ \ell_j\rangle.
\]
Expand
\[
\eta(y)\prod_{j=0}^n \zeta\circ\ell_j(y)
=
\int_{\xi\in\reals^m} a(\xi)e^{iy \cdot \xi}\,d\xi
\]
where $a(\xi) = O((1+|\xi|)^{-m-1})$.
It suffices to prove that
\[
\big|
\langle \scriptk,\ e^{iy\cdot\xi}\prod_j f_j\circ\ell_j(y)\rangle
\big|
\le AC\prod_j\|f_j\|_{p_j}.
\]
But since the intersection of the nullspaces of $\{\ell_j\}$ is trivial,
there exist linear mappings $L_j:\reals^\kappa\to\reals^m$ 
such that $y\equiv \sum_{j=0}^n L_j\circ\ell_j(y)$ for all $y\in\reals^m$.
Consequently $e^{iy\cdot\xi}$ may be rewritten as $\prod_j e^{i\ell_j(y)\cdot L_j^*(\xi)}$,
and if we define $\tilde f_j(x) = f(x)e^{ix\cdot L_j^*(\xi)}$ then our functional becomes
$\langle\scriptk,\ \prod_j\tilde f_j\circ\ell_j\rangle$.
Since $\|\tilde f_j\|_{p_j}=\|f_j\|_{p_j}$, \eqref{lemmaconclusion}
follows from \eqref{lemmahypothesis}.
\end{proof}

\begin{proof}[Proof of Theorem~\ref{thm:conditional}]
Let subspaces $\{V_j\}$, a degree $D\ge 1$, and a polynomial $P$ of degree $\le D$ be given.
We will prove the result only for $q\ge 1$, which permits the use of duality
and thus simplifies notation, leaving the simple modifications for $q<1$ to the reader.
Let $p_0$ be the exponent dual to $q$, and recall that $\ell_0(y)\equiv 0$.   
By duality, matters reduce to a multilinear functional
\[
\scriptt_P(f_0,\cdots,f_n) = \Big\langle K(y)e^{iP(x,y)},\,
\prod_{j=0}^n f_j(x+\ell_j(y))\Big\rangle,
\]
where the pairing is with respect to $(x,y)$.
We proceed by induction on $D$, the result for $D=0$ being given as a hypothesis.

Decompose $P=P_D+R$ where $R$ has degree $<D$ and $P_D$ is homogeneous of degree $D$.
If $P_D$ is degenerate then by decomposing $P_D = \sum_j q_j(x+\ell_j(y))$
for certain real-valued polynomials $q_j$, we may rewrite
$\scriptt_P(f_0,\cdots,f_n) = \scriptt_R(e^{iq_0}f_0,\cdots,e^{iq_n} f_n)$.
By induction on the degree, $\scriptt_R$ satisfies the desired estimate,
which is equivalent to the desired estimate for $\scriptt_P$ since
$\|e^{iq_j}f_j\|_{p_j}=\|f_j\|_{p_j}$.

If $P_D$ is nondegenerate, there exists a smallest integer $N$ such that
$P_D(2^{N}x,2^{N}y)$ has norm $\ge 1$ in the quotient space of
homogeneous polynomials of degree $D$ modulo degenerate homogeneous
polynomials of that degree, with respect to some fixed choice of norm on that 
finite-dimensional space.
This norm is then also $\le 2^{D}$.

By rescaling the variables $x,y$ by a factor of $2^N$
we may reduce henceforth to the case $N=0$.
Because $K$ is homogeneous of the critical degree and $\sum_{j=0}^n p_j^{-1}=1$,
such a rescaling does not affect the inequality to be proved.

By replacing $P_D$ by $P_D-\sum_j q_j(x+\ell_j(y))$ for appropriate
homogeneous real-valued polynomials $q_j$ of degree $D$,
we may assume henceforth that $P_D$ has norm $\lesssim 1$
in the space of all homogeneous polynomials of degree $D$, rather than
merely in the quotient space.

Fix a cutoff function $\zeta\in C^\infty_0(\reals^k)$
which is $\equiv 1$ in some neighborhood of the origin,
and consider the truncated multilinear functional
$\scriptt'_P = \langle \zeta(y) Ke^{iP(x,y)},
\prod_{j=0}^n f_j(x+\ell_j(y))\rangle$.
Introduce also $\eta\in C^\infty_0(\reals^d)$ such that
$\sum_{\nu\in \integers^{d}} \eta(x-\nu)\equiv 1$ 
and decompose
\[\scriptt'_P(f_0,\cdots,f_n)
= \sum_{\nu=(\nu_0,\cdots,\nu_n)\in\integers^{(n+1)d}}
\big\langle \zeta(y) Ke^{iP(x,y)},\  
\prod_{j=0}^n f_j(x+\ell_j(y))\eta(x+\ell_j(y)-\nu_j)\big\rangle.\]
Because of the presence of the compactly supported factor $\zeta(y)$,
there exists $C_0<\infty$ such that any term with parameter $\nu$
vanishes unless $\max_{j}|\nu_j-\nu_0|\le C_0$.
Therefore by the triangle inequality, H\"older's inequality
and the hypothesis $\sum_j p_j^{-1}=1$, it suffices to prove the desired bound
for fixed $\nu$, with $f_j$ replaced by $f_j\eta(x+\ell_j(y)-\nu_j)$,
so long as a majorization uniform in $\nu$ is obtained.

Make the change of variables $x=z+\nu_0$. 
$P(x,y) = P_D(x,y)+R(x,y)
= P_D(z,y) + \tilde R_{\nu_0}(z,y)$ where $\tilde R_{\nu_0}$ 
has degree $<D$. Moreover, if another cutoff function $\tilde\zeta\in C^\infty_0$
is chosen to be $\equiv 1$ in a sufficiently large neighborhood
of the origin, then our functional may be written as
\[
\big\langle
K(y)e^{i\tilde R_{\nu_0}(z,y)}
\big(e^{iP_D(z,y)}\zeta(y)\tilde\zeta(z,y)\big),
\  
\prod_{j=0}^n \tilde f_j(z+\ell_j(y))
\big\rangle
\]
where each $\tilde f_j$ is an appropriate translate of the product of $f_j$
with $\eta(x+\ell_j(y)-\nu_j)$.
Now the factor $\zeta(y)e^{iP_D(z,y)}\tilde\zeta(z,y)$ is smooth and compactly 
supported as a function of $(z,y)$, and is bounded above in any $C^s$ norm 
by a constant depending only on $s,D$
and the choice of $\tilde\zeta$. By induction on $D$, the functional
$\langle K e^{i\tilde R_{\nu_0}(z,y)},\ 
\prod_{j=0}^n \tilde f_j(z+\ell_j(y)) \rangle$
maps $\otimes_j L^{p_j}$ to $\complex$ boundedly,
uniformly for all polynomials $\tilde R_{\nu_0}$ of degree $<D$.
Therefore we may invoke Lemma~\ref{lemma:truncation} to conclude that
\[
\Big|
\big
\langle
Ke^{i\tilde R_{\nu_0}(z,y)}\big(e^{iP_D(z,y)}\zeta(y)\tilde\zeta(z,y)\big),
\  
\prod_{j=0}^n \tilde f_j(z+\ell_j(y))
\big\rangle
\Big|
\le C\prod_j\|f_j\|_{p_j},
\]
uniformly in $\nu$.

It remains to analyze 
\[
\langle Ke^{iP}(1-\zeta(y)),\ \prod_j f_j(x+\ell_j(y))\rangle
=\sum_{r=0}^\infty 
\langle Ke^{iP}h(2^{-r}y),\ \prod_j f_j(x+\ell_j(y))\rangle
\]
where $h(y) = \zeta(y/2)-\zeta(y)$.
We claim that the $r$-th summand is majorized by $C2^{-\eps r}\prod_j\|f\|_{p_j}$,
for some $\eps>0$  and $C<\infty$ depending only on $D,\{V_j\}$.
Because of the homogeneity of $K$ and the condition $\sum_j p_j^{-1}=1$, this is equivalent to
\begin{equation} \label{lastthing}
\big|\big\langle
K(y)h(y)e^{iP_r},\ \prod_j f_j(x+\ell_j(y))
\big\rangle\big|
\le C2^{-\eps r}\prod_j\|f\|_{p_j},
\end{equation}
via the substitution
$(x,y) = (2^r x',2^ry')$,
where $P_r(x,y) = P(2^rx,2^ry) = 2^{Dr}P_D(x,y)+ \tilde R_r(x,y)$
where $\tilde R$ has degree $<D$ and is real-valued.

As above, we may reduce matters to the case where all $f_j$ are supported in
a fixed bounded set, by introducing a partition of unity.
$\tilde R_r$ is thereby further modified and depends also on the index $\nu_0$; it
satisfies no useful upper or lower bounds but still has degree strictly $<D$.

Any such polynomial $2^{Dr} P_D + \tilde R$, where $\tilde R$ has degree $<D$,
has norm $\ge c2^{Dr}$ in the quotient space
of polynomials of degree $\le D$ modulo degenerate polynomials of degreee $\le D$,
for some fixed constant $c>0$, depending only on $D$.
The function $Kh$ is Lipschitz, since $K$ is a Calder\'on-Zygmund kernel
and $h$ vanishes identically in some neighborhood of the origin.
Therefore \eqref{lastthing} follows 
from the assumption that $\{V_j\}$ has the uniform decay property.
\end{proof}

\section{Sublevel sets}
\label{section:sublevel}

We next prove 
Theorem~\ref{thm:sublevel}. The argument requires an elaboration of
the characterization of nondegeneracy in terms of difference operators
that was established in Lemma~\ref{lemma:differences}. 

Define divided difference operators
\begin{align*}
&\delta_{j,r}f(x) = [f(x+re_j)-f(x)]/r,
\\
&\Delta_{\alpha,r} = 
\delta_{1,r_{1,1}}\circ\delta_{1,r_{1,2}}\cdots
\circ\delta_{1,r_{1,\alpha_1}}
\circ\cdots\circ\delta_{n,r_{n,1}}\circ\cdots
\delta_{n,r_{n,\alpha_n}}
\end{align*}
where in the second definition $r$ has arbitrary components $r_{j,k}\in\reals^+$ 
for all $1\le j\le m$ and $1\le k\le\alpha_j$.

Let subspaces $V_j\subset\reals^m$ of some dimension $\kappa$ be given for $1\le j\le n$,
where $n,m,\kappa$ are arbitrary. Let
$P$ be any polynomial which is nondegenerate relative to $\{V_j\}$.
Let $D$ be the degree of $P$.
As already shown, we may suppose without loss of generality that 
the homogeneous component of $P$ of degree $D$ is itself nondegenerate.

Let $L = \sum_{|\beta|=D} b_\beta \Delta_\beta$ be the homogeneous difference operator
constructed in Lemma~\ref{lemma:differences}, which annihilates all degenerate
polynomials $p_j \circ\pi_j$, but does not annihilates $P$.
Consider the more general operators
\begin{equation} \label{moregeneralL}
L_r = \sum_{|\beta|=D} b_\beta \prod_{j,k} \Delta_{\beta,r_\beta}
\end{equation}
where there is associated to each index $\beta$ a separate multi-parameter
$r_\beta$ whose components are $r_{\beta,j,k}\in\reals^+$.
The proof of Lemma~\ref{lemma:differences} demonstrates that
\begin{equation}   \label{moreflexible}
L_r (P-\sum_j f_j\circ\pi_j)\equiv 1  \ \text{a.e.\ for all $f_j\in L^1_{\text{loc}}$},
\end{equation}
for all multi-parameters $r$. 
(We require that $f_j$ be locally integrable, rather than merely measurable,
in order that it can be interpreted as a distribution.)
Indeed, for any homogeneous polynomial $Q$ of degree $D$,
$L_r Q(0)=LQ(0)$ for any $r$, by homogeneity because we are working
with {\em divided} difference operators.
On the other hand, $L_rQ(0)=0$ for any homogeneous polynomial
of any other degree. 
The rest of the proof parallels that of Lemma~\ref{lemma:differences}.

For any function $g$, point $x$, multi-index $\beta$ and 
multi-parameter $r_\beta$,
$\Delta_{\beta,r_\beta}g(x)$ is a linear combination of
$2^{|\beta|}$ values $g(x+\sum_{i} r_{\beta,i}\sigma_{\beta,i}y_{\beta,i})$,
where $i$ ranges over an index set of cardinality $|\beta|=D$,
each $y_{\beta,i}$ belongs to $\reals^m$,
and each $\sigma_{\beta,k}\in\{0,1\}$,
with one summand for each of the $2^D$ pairs $(i,\sigma_{\beta,i})$.
Moreover, the coefficients in this linear combination
are fixed constants, depending on $(\beta,i,\sigma_{\beta,i})$,
times $\prod_i r_{\beta,i}^{-1}$.
Each $y_{\beta,i}$ is in fact one of the $m$ unit coordinate vectors
in $\reals^m$, but this information will not be used in the further discussion.

$L_r$ involves a further summation over $\beta$, so that 
$L_rf(x)$ is a linear combination of values of $f$ at the points
of $x+\cup_\beta \{\sum_i r_{\beta,i}\sigma_{\beta,i}y_{\beta,i}\}$. 
It will be convenient for the reasoning below to simplify the form of this
linear combination by placing the variables $\beta,i$ on a more even footing.
To achieve this we enlarge the collection of points 
$\cup_\beta \{\sum_i r_{\beta,i}\sigma_{\beta,i}y_{\beta,i}\}$, 
whose cardinality
is $N+2^D$ where $N$ is the number of indices $\beta$,
by forming the finite set
\begin{equation}
Y_r = \big\{\sum_\beta\sum_i r_{\beta,i}\sigma_{\beta,i}y_{\beta,i}: 
\text{ each } \sigma_{\beta,i}\in\{0,1\}\big\},
\end{equation}
whose cardinality is $2^{N+D}$.
With this notation,
$L_r$ takes the form
\begin{equation} \label{revisedLsubr}
L_rg(x) = \sum_{y\in Y_r} c_{y,r} g(x+y)
\end{equation}
where each coefficient $c_{y,r}$ is proportional to 
the product of reciprocals of certain components of $r$,
while $Y_r$ takes the form
$Y_r = \{\sum_{\alpha\in\scripta} r_\alpha\sigma_\alpha y_\alpha: \sigma\in \{0,1\}^{|\scripta|}\}$
for a certain index set $\scripta$ and points $y_\alpha$.
For typical (non-monomial) $L_r$, 
the overwhelming majority of the coefficients $c_{y,r}$ will equal zero.

Given $P$ and functions $f_j$, let
\begin{equation}
E_\eps=\big\{x\in\reals^k: |P(x)-\sum_j f_j(\pi_j(x))|<\eps\big\}.
\end{equation}
In the spirit of \cite{ccw},
we now describe a certain combinatorial restriction on
these sublevel sets which is implied by nondegeneracy of $P$.
\begin{lemma} \label{notszemeredi}
Suppose that a homogeneous polynomial $P:\reals^m\to\reals$ is nondegenerate
relative to a collection of subspaces $\{V_j\}$.
Let the sets $Y_r$ be defined as above, where the difference
operators $L_r$ satisfy \eqref{moreflexible}. 
Then for any Lebesgue measurable functions $f_j$, 
any multi-parameter $r$,
and almost every $x\in\reals^m$, 
\begin{equation}
\text{If $x+y\in E_\eps$ for every $y\in Y_r$ then 
at least one component of $r$ is $\le C\eps^{1/D}$.}
\end{equation}
\end{lemma}

\begin{proof}
It suffices to prove this for locally integrable functions $f_j$,
since the general case then follows from a limiting argument.
Let $g = P-\sum_j f_j\circ\pi_j$.
If $x+y\in E_\eps$ for every $y\in Y_r$ then
each term in the sum \eqref{revisedLsubr} is $O(\rho^{-D}\eps)$,
where $\rho$ is the smallest component of $r$.
Thus $1 = L_r g(x)=O(\rho^{-D}\eps)$,
so $\rho\le C\eps^{1/D}$.
\end{proof}

We now abstract the situation to which the proof of Theorem~\ref{thm:sublevel}
has been reduced.
Let $Y\subset\reals^m$ be an arbitrary finite, nonempty subset.  Write
$Y=\{y_\alpha: \alpha\in\scripta\}$ where $\scripta$ is a finite index set.
To $Y$ associate the sets 
$Y_r = \{\sum_{\alpha\in\scripta} r_\alpha\sigma_\alpha y_\alpha: \sigma\in \{0,1\}^{|\scripta|}\}$,
where $r=(r_\alpha)$, $\sigma=(\sigma_\alpha)$,
each $r_\alpha\in\reals^+$, and each $\sigma_\alpha\in\{0,1\}$.
Thus $Y_r$ has cardinality $2^{|\scripta|}$.

\begin{lemma} \label{lemma:combinatorial}
Let $Y$ be any nonempty finite subset of $\reals^m$.
Let $E\subset\reals^m$ be measurable and contained in a fixed bounded set.
Let $\eps\in(0,1]$.
Suppose that for every $r\in(0,1]^{\scripta}$, for almost every $x\in\reals^m$,
if $x+y\in E$ for every $y\in Y_r$ then at least one
component $r_\alpha$ is $\le\eps$.
Then $|E|\le C\eps^\delta$
where $C,\delta\in\reals^+$, $\delta$ depends only
on the cardinality of $Y$, and $C$ depends only on $Y$. 
\end{lemma}

\begin{proof}
The special case where $Y_r$ is simply the set of all vectors
$\sum_{j=1}^m r_j\sigma_j e_j$, where $e_j$ are the coordinate vectors,
and $r=(r_1,\cdots,r_m)$, was treated in \cite{ccw}.
The general case may be reduced to that special case by the following lifting argument.

Introduce $\reals^M=\reals_x^m\times\reals_t^{|\scripta|}$, 
adding one real coordinate for each index $\alpha$. 
Let $e_\alpha\in\reals^{|\scripta|}$ be the unit vector corresponding to the $\alpha$-th
coordinate.
Define $E^\dagger=E\times B$ where $B$ is a fixed large
ball in $\reals^{|\scripta|}$.

Introduce the shear transformation $T:\reals^M\to\reals^M$
defined by $T(x,t) = (x-\sum_\alpha t_\alpha y_\alpha,t)$, 
where $t=\sum_\alpha t_\alpha e_\alpha$,
and let $E^\ddagger = T(E^\dagger)$.
Then for any $r,\sigma,s$  and almost every $x$,
\[
x+\sum_\alpha r_\alpha\sigma_\alpha y_\alpha \in E
\text{ if and only if } 
T(x,s)+(0,\sum_\alpha r_\alpha \sigma_\alpha e_\alpha)\in E^\ddagger.
\]
Indeed, 
$x+\sum_\alpha r_\alpha \sigma_\alpha y_\alpha\in E$ 
is equivalent to 
$(x+\sum_\alpha r_\alpha \sigma_\alpha y_\alpha,s+\sum_\alpha r_\alpha \sigma_\alpha e_\alpha)
\in E^\dagger$,
since $E^\dagger$ is invariant under translations in the second set of variables.
Next
\begin{align*}
T(x+\sum_\alpha r_\alpha \sigma_\alpha y_\alpha,
&s+\sum_\alpha r_\alpha \sigma_\alpha e_\alpha) 
\\
&= (x+\sum_\alpha r_\alpha \sigma_\alpha y_\alpha
- \sum_\alpha (s_\alpha + r_\alpha \sigma_\alpha )y_\alpha,
s+\sum_\alpha r_\alpha \sigma_\alpha e_\alpha) 
\\
&= (x-\sum_\alpha s_\alpha y_\alpha ,
s+\sum_\alpha r_\alpha \sigma_\alpha e_\alpha )
\\
&= T(x,s)+(0,\sum_\alpha r_\alpha \sigma_\alpha e_\alpha). 
\end{align*}
Thus 
\begin{align*}
x+\sum_\alpha r_\alpha \sigma_\alpha y_\alpha \in E
&\Leftrightarrow
(x+\sum_\alpha r_\alpha \sigma_\alpha y_\alpha,
s+\sum_\alpha r_\alpha \sigma_\alpha e_\alpha)\in E^\dagger
\\
&\Leftrightarrow T(x+\sum_\alpha r_\alpha\sigma_\alpha y_\alpha,
s+\sum_\alpha r_\alpha \sigma_\alpha e_\alpha )\in E^\ddagger
\\
&\Leftrightarrow T(x,s)+(0,\sum_\alpha r_\alpha \sigma_\alpha e_\alpha)\in E^\ddagger.
\end{align*}

Let $Y^\ddagger_r = \{(0,\sum_\alpha r_\alpha \sigma_\alpha e_\alpha)\}\subset\reals^M$.
Suppose now that whenever $x+y\in E$ for every $y\in Y_r$,
some component $r_\alpha$ is $\le\eps$. Then for any $z\in \reals^M$,
if $z+y\in E^\ddagger$ for every $y\in Y^\ddagger_r$,
then some $r_\alpha\le\eps$.

By Fubini's theorem, we may freeze the variable $x$, reducing matters
to $\reals^{|\scripta|}$; we are given (for almost every $x\in\reals^m$)
a set $E^*$, contained in a fixed bounded subset of $\reals^{|\scripta|}$, such that 
for any $z\in \reals^{|\scripta|}$,
if $z+\sum_\alpha r_\alpha \sigma_\alpha e_\alpha \in E^*$ for every $\sigma\in\{0,1\}^{|\scripta|}$
then some $r_\alpha\le\eps$.
It was shown in \cite{ccw} that this forces $|E^*|\le C\eps^\delta$
for some $\delta>0$. See Lemma 3.7 of \cite{gowerscurrent} for a more precise 
discrete analogue which implies the continuum version.
\end{proof}

\begin{proof}[Proof of Theorem~\ref{thm:sublevel}]
Let $P,\{V_j\},g_j,\eps$ be as in Theorem~\ref{thm:sublevel},
and let $D$ be the degree of $P$.

In proving the theorem,
we may suppose without loss of generality that the homogeneous
part of $P$ of degree $D$ is nondegenerate with respect to $\{V_j\}$,
for if not, it may be removed by replacing $P$ by $P-\sum_j p_j\circ\pi_j$
for suitable polynomials $p_j$ of degree $D$.
Let the difference operators $L_r$ and associated finite
sets $Y$ be associated to $P_D$
as in Lemma~\ref{notszemeredi} and the accompanying discussion.

Set $F = P-\sum_j g_j\circ\pi_j$ and consider its sublevel set
$E_\eps=\{z: |F(z)|<\eps\}$.
The difference operators $L_r$ used in the proof of Lemma~\ref{notszemeredi}
satisfy $L_r(P_D-\sum_j g_j\circ\pi_j)\equiv 1$,
but they annihilate all polynomials of degrees $<D$, so 
also satisfy $L_r(P-\sum_j g_j\circ\pi_j))\equiv 1$.
Therefore the proof of Lemma~\ref{notszemeredi} demonstrates
that (disregarding a set of measure zero)
if $x+y\in E_\eps$ for all $y\in Y_r$ then
the smallest component of $r$ is $\lesssim \eps^{1/D}$.
It now suffices to invoke Lemma~\ref{lemma:combinatorial}.
\end{proof}

\section{Further results}
Theorem~\ref{thm:lowcodimension} is a direct consequence of Lemma~\ref{lemma:lowcodimension} together with
Theorem~\ref{thm:simpledecay}.
The following extension of Theorem~\ref{thm:codimone} will be used to derive 
Theorem~\ref{thm:noslicing}.

\begin{theorem} \label{thm:augment}
Let $\{V_j\}$ be any finite collection of subspaces of $\reals^m$,
and let $\{W_i\}$ be any finite collection of codimension one subspaces of $\reals^m$.
If $\{V_j\}$ has the uniform decay property, then so does $\{V_j\}\cup\{W_i\}$.
\end{theorem}

\begin{proof}
It suffices to prove this in the case where a single codimension one subspace
$W$ is given. We will first prove that any nondegenerate polynomial
has the decay property, then address the uniformity issue at the end of the proof.
Choose coordinates with respect to which $W=\{(x',x_m): x_m=0\}$.
Let a polynomial $P$ be nondegenerate relative to the augmented
collection $\{V_j\}\cup\{W\}$.
We may assume that no subspaces $V_{j}$ are contained in $W$, for any such
spaces may be deleted from $\{V_j\}$ without affecting the nondegeneracy of $P$.

$\p P/\p x_m$ is nondegenerate relative to $\{V_j\}$. For if not,
then there exists a polynomial decomposition $\p P/\p x_m= \sum_j q_j\circ\pi_j$.
Since no space $V_j$ is contained in $W$, $\p(Q_j\circ\pi_j)/\p x_m
= (v_j\dot\nabla Q_j)\circ\pi_j$ for certain nonzero vectors $v_j$.
Therefore there exist polynomials $Q_j$ such that $\p (Q_j\circ\pi_j)/\p x_m
= q_j\circ\pi_j$, and hence by setting $\tilde P = \sum_j Q_j\circ\pi_j$
we have $\p (P-\tilde P)/\p x_m \equiv 0$. Thus $P - \sum_j Q_j\circ\pi_j$
is a function of $x_m$ alone, whence $P$ is degenerate.

Let $d$ be the degree of $P$.
It now follows that there exists $z\in\reals^m$ for which
$P_z(x) = P(x',x_m)-P(x',x_m+z)$ is nondegenerate relative to $\{V_j\}$.
If we consider the quotient space $\scriptp$ of all polynomials of
of degrees $\le d$ modulo the subspace of all such polynomials
which are degenerate relative to $\{V_j\}$ with an inner product structure, then
$\|P_z\|_{\scriptp}^2$ is a polynomial in $z$ which does not vanish identically.
Hence there exist $C,\delta\in\reals^+$
such that for any ball $B$ of fixed finite radius, for any $\eps>0$,
\begin{equation} \label{augsublevel}
|\{z\in B: \|P_z\|_{\scriptp}^2<\eps\}|\le C\eps^\delta.
\end{equation}
We may now argue as in the proof of Theorem~\ref{thm:slices} to conclude the proof
that $P$ has the decay property.

The same reasoning as above demonstrates that if a family $\{P_\alpha\}$ of polynomials
of uniformly bounded degrees is uniformly nondegenerate relative to $\{V_j\}$,
then $\{\p P_\alpha/\p x_m\}$ is uniformly nondegenerate relative
to $\{V_j\}\cup\{W\}$. Hence $\|\p P_\alpha/\p x_m\|_{\scriptp}$ is bounded
away from zero, uniformly in $\alpha$.
From this the sublevel set bound \eqref{augsublevel} follows by elementary reasoning.
\end{proof}

\begin{corollary}
Let a real-valued polynomial $P$ be nondegenerate relative
to a finite collection $\{V_j\}$ of subspaces of $\reals^m$,
and suppose that $V_j$ has codimension one for all $j>1$.
Then $\{V_j\}$ has the uniform decay property. Moreover
\begin{equation}  \label{Ltwoversion}
|\Lambda_\lambda(f_1,\cdots,f_n)|\le C|\lambda|^{-\eps}\|f_1\|_2 \prod_{j>1}\|f_j\|_\infty,
\end{equation}
with uniform bounds if $P$ belongs to a family of uniformly nondegenerate polynomials.
\end{corollary}

Indeed, this follows from the $L^\infty$ bound  by interpolating
with the trivial bound $\lesssim \|f_1\|_1\prod_{j>1}\|f_j\|_\infty$.

\begin{proof}[Proof of Theorem~\ref{thm:noslicing}]
We will prove power decay in the stronger form \eqref{Ltwoversion}.
The case $n=1$ is a well-known fact, as discussed earlier.
Let $W$ be the span of those $V_j$ with $2\le j\le M+1$,
and $\tilde W$ be the span of those with $M+1<j\le n$,
$n\le 1+2M$ and $(1+M)\kappa\le m$.
By the general position hypothesis, $W$ has dimension $M\kappa$
and $\tilde W$ has dimension $(n-M-1)\kappa$; both of these dimensions
are $\le m-1$.

Let functions $f_j\in \lt(V_j)$ be given.
As usual, we may assume each $f_j$ to be supported in a fixed bounded subset of $V_j$.

The proof is divided into two main cases, depending on whether or not $f_1$
is $\lambda$-uniform, as defined in Definition~\ref{defn:uniformity}.
Let $A(\lambda)$ be the best constant in the inequality \eqref{Ltwoversion}.
If $f_1$ is not $\lambda$-uniform, then it follows by induction on $n$,
as in the proof of Theorem~\ref{thm:slices}, that
$A(\lambda) \le C|\lambda|^{-\delta} + (1-|\lambda|^{-\tau})A(\lambda)$.

Define functions $F,\tilde F$ on $W,\tilde W$ respectively by the relations
$F\circ\pi = \prod_{j=2}^{M+1} f_j\circ\pi_j$,
$\tilde F\circ\tilde\pi = \prod_{j>M+1} f_j\circ\pi_j$
where $\pi,\tilde\pi$ denote the orthogonal projections from $\reals^m$ onto $W,\tilde W$
respectively.
Then by general position,
$\|F\|_2\lesssim\prod_{j=2}^{M+1}\|f_j\|_2$
and
$\|\tilde F\|_2\lesssim\prod_{j>M+1}\|f_j\|_2$.

Consider next the case where $f_1$ is $\lambda$-uniform.
In this case the hypothesis that $P$ is nondegenerate plays no role.
Our integral is
\[
\Lambda_\lambda=
\int_{\reals^m} e^{i\lambda P(x)}(F\circ\pi)(\tilde F\circ\pi)(f_1\circ\pi_1) \eta.
\]
Now $P$ may be nondegenerate relative to the collection of three subspaces $V_1,W,\tilde W$.
If it is, then the proof is complete by virtue of the preceding corollary,
in its more precise version \eqref{Ltwoversion}.
Thus we may assume henceforth that $P$  is decomposable
as $P = p\circ\pi + \tilde p\circ\tilde\pi+q\circ\pi_1$ for certain polynomials $p,\tilde p,q$.

Set  
$G = e^{i\lambda p}F$,
$\tilde G = e^{i\lambda \tilde p}\tilde F$,
and $g = e^{i\lambda q}f_1$. Our integral becomes
\[
\int (G\circ\pi)(\tilde G\circ\tilde\pi)(g\circ\pi_1)\eta
= c\iiint \widehat{G}(\xi)\widehat{\tilde G}(\tilde\xi)\widehat{g}(\zeta)
\widehat{\eta}(\pi^*\xi+\tilde\pi^*\tilde\xi+\pi_1^*\zeta)\,d\xi\,d\tilde\xi\,d\zeta
\]
where $\pi_1^*$ is the adjoint of $\pi_1:\reals^m\to V_1$, and analogously for $\pi^*,\tilde\pi^*$.
The $\lambda$-uniformity condition gives a bound $\|\widehat{g}\|_{L^\infty}=O(|\lambda|^{-\tau})$
for some fixed $\tau>0$. Therefore for $|\lambda|\ge 1$
\[ 
|\Lambda_\lambda|
\lesssim |\lambda|^{-\tau}\|f_1\|_2
\iiint |\widehat{G}(\xi)\widehat{\tilde G}(\tilde\xi)|
\cdot
|\widehat{\eta}(\pi^*\xi+\tilde\pi^*\tilde\xi+\pi_1^*\zeta)|\,d\xi\,d\tilde\xi\,d\zeta.
\]
It therefore suffices to have 
\begin{equation}  \label{etahatdecay}
\sup_{\tilde\xi} \int 
|\widehat{\eta}(\pi^*\xi+\tilde\pi^*\tilde\xi+\pi_1^*\zeta)|\,d\xi\,d\zeta<\infty
\end{equation}
and the same with the roles of $\tilde\xi,\xi$ interchanged.
The mapping $\reals^{\kappa+M\kappa}\owns (\zeta,\xi)\to \pi^*\xi+\pi_1^*\zeta\in\reals^m$
is linear and injective by the general position hypothesis, since $(1+M)\kappa\le m$,
and $\eta$ may be taken to belong to $C^{K}$ for any preassigned $K$, so that $\widehat{\eta}$
decays rapidly. Thus \eqref{etahatdecay} holds.
Since the roles  of $\xi,\tilde\xi$ are symmetric,
it continues to hold when they  are interchanged.
\end{proof}


More general results, in which the dimensions of the subspaces $V_j$ are not
all required  to be equal, can be proved in the same way.


\begin{thebibliography}{9}

\bibitem{ccw} 
A.~Carbery, M.~Christ, and J.~Wright,
{\em Multidimensional van der Corput and sublevel set estimates},
J. Amer. Math. Soc. 12 (1999), no. 4, 981--1015.

\bibitem{gowers} 
W.~T.~Gowers,
{\em  A new proof of Szemer\'edi's theorem for arithmetic progressions of length four},
Geom. Funct. Anal. 8 (1998), no. 3, 529--551.

\bibitem{gowerscurrent}
\bysame,
{\em Arithmetic progressions in sparse sets},  
Current developments in mathematics, 2000,  149--196, Int. Press, Somerville, MA, 2001.


\bibitem{bilinearHilbert}
M.~Lacey and C.~Thiele, {\em $L^p$ estimate on the bilinear Hilbert 
transform}, Ann. of Math. 146 (1997), 693-724.

\bibitem{bHilbert2}
\bysame,
{\em On Calder\'on's conjecture}, Ann. of Math. 149 (1999), 475-496.

\bibitem{mtt}
C.~Muscalu, T.~Tao, and C.~Thiele, 
{\em Multi-linear operators given by singular multipliers },
J. Amer. Math. Soc. 15 (2002), 469-496.

\bibitem{phongstein}
D.~H.~Phong and E.~M.~Stein, {\em 
Operator versions of the van der Corput lemma and Fourier integral operators},
Math. Res. Lett. 1 (1994),  27-33.

\bibitem{phongsteinsturm}
D.~H.~Phong, E.~M.~Stein, and J.~Sturm,
{\em Multilinear level set operators, oscillatory integral operators, and Newton polyhedra},
Math. Ann. 319 (2001), no. 3, 573--596.

\bibitem{riccistein}
F.~Ricci and E.~M.~Stein, {\em Harmonic analysis on nilpotent groups and 
singular integrals }, J. Funct. Anal. 73 (1987), 179-194.

\bibitem{steinbook}
E.~M.~Stein, 
{\em Harmonic analysis: real-variable methods, orthogonality, and oscillatory integrals. With the assistance of Timothy S. Murphy}, 
Princeton Mathematical Series, 43. Monographs in Harmonic Analysis, III. 
Princeton University Press, Princeton, NJ, 1993. 


\end{thebibliography}
\end{document}